\documentclass[11pt]{article}
\usepackage{graphics}
\usepackage{amssymb}
\usepackage{pstricks}
\usepackage[french]{babel}
\newcommand{\ddate}{31 octobre 2002}

\newtheorem{thm}{Th\'eor\`eme}[section]

\newtheorem{ex}[thm]{Example}
\newtheorem{Theorem}[thm]{Th\'eor\`eme}

\newtheorem{Lemme}[thm]{Lemme}
\newtheorem{Proposition}[thm]{Proposition}
\newtheorem{Remarque}[thm]{Remarque}

\newtheorem{Corollaire}[thm]{Corollaire}

\newtheorem{ccote}[thm]{}

\newcommand{\preu}{\noindent {\sc Preuve: \ }}

\newcommand{\cqfd}{\unskip\kern 6pt\penalty 500
\raise -2pt\hbox{\vrule\vbox to10pt{\hrule width
 4pt\vfill\hrule}\vrule}\smallskip}

\newcommand{\bbr}{{\mathbb{R}}}
\newcommand{\bbc}{{\mathbb{C}}}

\newcommand{\bbn}{{\mathbb{N}}}

\newcommand{\calc}{{\cal C}}
\newcommand{\cald}{{\cal D}}

\newcommand{\calh}{{\cal H}}

\newcommand{\call}{{\cal L}}

\newcommand{\caln}{{\cal N}}
\newcommand{\calo}{{\cal O}}

\newcommand{\calt}{{\cal T}}

\newcommand{\pcirc}{\kern .7pt {\scriptstyle \circ} \kern 1pt}
\newcommand{\iso}{\cong}

\newcommand{\eqref}[1]{(\ref{#1})}
\newcounter{exo}

\newcommand{\mun}{{-1}}
\newcommand{\grad}{{\rm grad}}
\newcommand{\sk}[1]{\vskip #1 mm}

\newcommand{\llangle}[2]{\langle #1 ,#2 \rangle}
\newcommand{\vect}{{\rm Vec\,}}
\newcommand{\nua}[2]{\caln^{#1}_{#2}}
\newcommand{\lacets}[1]{{\rm Lac\,}_{#1}}
\newcommand{\lie}{{\rm Lie\,}}
\newcommand{\moeb}[1]{{\rm M\ddot{o}b}_{#1}}
\newsavebox{\trait}
\newsavebox{\brasa}
\newsavebox{\brasb}
\newsavebox{\conab}
\newsavebox{\conac}
\newsavebox{\conad}
\newsavebox{\conbc}
\newsavebox{\conbd}
\newsavebox{\concd}
\newsavebox{\pbras}

\title{Contr\^ole des bras articul\'es et \\
transformations de M\"obius}
\author{Jean-Claude HAUSMANN}
\date{\ddate}

\begin{document}
\maketitle
\begin{abstract}
Pour un $m$-uple $a=(a_1,\dots,a_m)$ de nombres r\'eels positifs,
le bras articul\'e de type $a$ dans $\bbr^d$ est l'application
$f^a:(S^{d-1})^m\to\bbr^d$ d\'efinie par
$f^a(z_1,\dots,z_m)=\sum_{j=1}^ma_jz_j$.
On se propose d'attaquer le probl\`eme inverse via les relev\'es
horizontaux pour la distribution $\Delta^a$ orthogonale aux fibres de $f^a$.
On montre que les composantes connexes par courbes horizontales
sont les orbites d'une action sur $(S^{d-1})^m$ d'un produit de groupes
de transformations de M\"obius. On d\'etermine \'egalement, dans plusieurs cas,
les orbites d'holonomie de la distribution $\Delta^a$.
\end{abstract}
%\tableofcontents
\iffalse
\vskip 3cm
\setlength{\unitlength}{.05mm}%{1mm}%
\begin{pspicture}(0,0)(0,0)
%\pspicture[0](0,0)
%\color{gristq}\graphpaper[2](0,-10)(180,100)\color{black}
%\put(-8,0){\line(1,0){16}}
%\multiput(-9,0)(3,0){7}{\line(1,-1){3}}
\end{pspicture}
\end{document}
\fi

\section{Introduction} \label{intro}

Soit $f: M\to N$ une application diff\'erentiable ($\calc^\infty$)
entre deux vari\'et\'es $M$ et $N$. Le {\it probl\`eme inverse} pour $f$
consiste \`a trouver un proc\'ed\'e pour relever dans $M$ les courbes lisses
par morceau dans $N$. Plus pr\'ecis\'ement, \'etant donn\'e une
courbe lisse par morceau $c(t)\in N$ et un point
$q^0\in f^\mun(c(0))$, on aimerait disposer canoniquement
d'une courbe lisse par morceau $\tilde c(t)\in M$ telle que
 $\tilde c(0)=q^0$ et $f(\tilde c(t))=c(t)$.

Si $M$ est une vari\'et\'e riemannienne et $f$ est une
submersion propre, on peut alors proc\'eder
comme Ehresmann \cite[\S\,1]{Eh}. On consid\`ere une
{\it distribution} $\Delta^f$,
c'est-\`a-dire un champ de sous-espace vectoriels
$\Delta^f_q\subset T_qM$ pour tout $q\in M$,
en d\'efinissant $\Delta^f_q$ comme le compl\'ement orthogonal \`a
${\rm ker\,}T_qf$, o\`u
$T_qf:T_qM\to T_{f(q)}N$ est l'application tangente \`a $f$
en $q$. Une courbe lisse par morceau $b(t)\in M$
dont le vecteur vitesse $\dot b$
satisfait \`a $\dot b(t)\in\Delta^f_{b(t)}$ est appel\'ee
une {\it $\Delta^f$-courbe}, ou {\it courbe horizontale}.
Le probl\`eme inverse pour $f$ est alors r\'esolu de
la mani\`ere suivante~:
une courbe $c(t)\in N$ d\'efinie sur un intervalle
au voisinage de $0$ admet un unique rel\`evement
horizontal $\tilde c(t)\in M$ partant d'un point
$\tilde c(0)\in f^\mun(c(0))$.
En effet, la d\'emonstration est
la m\^eme que pour un fibr\'e diff\'erentiable muni d'une connexion
(voir \cite[lemme 2 p. 212]{DNF}). Rappelons que c'est
le principe de la preuve du th\'eor\`eme d'Ehresmann qui
dit qu'une submersion propre est un fibr\'e diff\'erentiable
\cite[\S\,1]{Eh}.

Pour une application diff\'erentiable $f:M\to N$ quelconque,
avec $M$ riemannienne, on peut tout de m\^eme consid\'erer
la distribution $\Delta^f$ où $\Delta^f_q$ est
le complément orthogonal à ${\rm ker\,}T_qf$. L'unicit\'e du
rel\`evement horizontal d'une courbe $c(t)\in N$ est alors
pr\'eserv\'ee puisque $T_qf$ envoie $\Delta^f_q$ isomorphiquement
sur l'image de $T_qf$. Son existence est garantie en tout cas si
$c(t)$ reste dans les valeurs r\'eguli\`eres de $f$, au dessus
desquelles $f$ est une submersion.
Cette solution horizontale du probl\`eme inverse est celle qui minimise
en chaque instant le co\^ut infinit\'esimal~:
$\tilde c(t)$ est horizontal si et seulement si
$\dot{\tilde c}(t)$ est le vecteur au dessus de
$\dot{c}(t)$ de norme minimale.
Observons que $\Delta^f$ est
une distribution de dimension non-constante~:
${\rm dim\,}\Delta^f_q={\rm rang\,}(T_qf)$.
Par exemple, si $N=\bbr^p$
avec la m\'etrique standard et que l'on \'ecrit
$f(q):=(x_1(q),\dots ,x_p(q))$, alors $\Delta^f$ est la distribution
engendr\'ee par les champs de gradients $\{\grad\, x_i\}_{i=1,\dots ,p}$.
Plus g\'en\'eralement, pour $f:M\to N$, la distribution $\Delta^f$
est engendr\'ee par les champs de gradient des fonctions
$\phi\pcirc f$ pour toutes les fonctions $\phi:N\to \bbr$.

Comme distribution sur une vari\'et\'e riemannienne,
$\Delta^f$ d\'efinit sur $M$ une
{\it g\'eom\'etrie sous-riemannienne} (\cite{Be},\cite{Mo}) ou {\it g\'eom\'etrie
de Carnot-Carath\'eodory} \cite{Gr}. En tant que distribution
engendr\'ee par des champs de vecteurs,
$\Delta^f$ d\'efinit \'egalement sur $M$ un
probl\`eme de {\it th\'eorie du contr\^ole} \cite{Ju}.
Ces th\'eories s'accordent sur l'utilit\'e du th\'eor\`eme de Sussmann
\cite[Th.\,4.1]{Su} qui d\'ecrit, pour une distribution
sur $M$, les composantes connexes par $\Delta$-courbes comme les
feuilles d'une distribution int\'egrable associ\'ee \`a $\Delta$.
Il s'agit encore d'un feuilletage dont les feuilles ne sont
pas toutes de m\^eme dimension, l'exemple standard \'etant
donn\'e par les orbites d'une action
diff\'erentiable d'un groupe de Lie sur $M$.
Dans le cas qui nous occupe d'une distribution $\Delta^f$,
obtenir une param\'etrisation de ces feuilles peut \^etre un
pas pr\'eliminaire important pour la r\'esolution du probl\`eme
inverse car il r\'eduit le nombre
de variables dans les \'equations diff\'erentielles.

Le but de cet article est
d'\'etudier les questions ci-dessus dans le cas o\`u
$M=(S^{d-1})^m$ est le produit de $m$ sph\`eres standard $S^{d-1}$
de rayon 1 (muni de la m\'etrique produit) et
$f:(S^{d-1})^m\to\bbr^d$ est un bras articul\'e.
Si $a:=(a_1,\dots ,a_m)\in\bbr^m_{\geq 0}$,
Le {\it bras articul\'e} (de type $a$) dans $\bbr^d$ est l'application
$f^a: (S^{d-1})^m\to\bbr^d$ d\'efinie par
$$f^a(z_1,\dots ,z_m):=\sum_{j=1}^m a_jz_j.$$
Les points de $(S^{d-1})^m$ seront appel\'es des {\it configurations}.
On peut s'imaginer un bras de robot dans $\bbr^d$, partant de l'origine,
form\'e de $m$ segments de longueurs $a_1,\dots,a_m$ et dont l'extr\'emit\'e
est $f^a(z)$.
Notons que certaines longueurs $a_i$ peuvent \^etre nulles,
ce qui nous sera techniquement utile dans certaines d\'emonstrations.

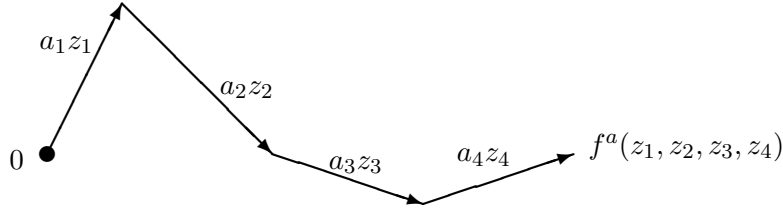
\setlength{\unitlength}{1mm}
\begin{picture}(60,40)(0,4)
\thicklines
\put(10,20){\vector(1,2){10}}
\put(20,40){\vector(1,-1){20}}
\put(40,20){\vector(3,-1){20}}
\put(60,13.3){\vector(3,1){20}}
\put(10,20){\circle*{2}}
\put(5,18){$0$}
\put(9,34){$a_1z_1$}
\put(33,28){$a_2z_2$}
\put(47.5,18){$a_3z_3$}
\put(64.5,19){$a_4z_4$}
\put(82,20){$f^a(z_1,z_2,z_3,z_4)$}
\end{picture}

\sk{-3}
Dans le cas d'un bras articul\'e,
le probl\`eme inverse prend la signification qu'il
a habituellement en robotique~:
trouver, \`a partir d'une
configuration initiale $\tilde c(0)$ du bras, un mouvement
$\tilde c(t)$
des diff\'erents segments faisant que $f^a(\tilde c(t))=c(t)$.

Le lecteur trouvera une \'etude des bras articul\'es dans \cite{Ha},
par exemple la caract\'erisation des points critiques de $f^a$
qui sont les configurations align\'ees $z_i=\pm z_j\ \forall\, i,j$
(voir aussi la remarque \ref{vcrit} plus loin),
ainsi qu'une premi\`ere classification des
pr\'eimages $M_b:=(f^a)^\mun(\{b\})$ pour une valeur r\'eguli\`ere $b\in\bbr^d$.
Les espaces $M_b$ sont reli\'es aux espaces de polygones. Par exemple,
si $d=2$ et $b\not= 0$, l'espace $M_b$ s'identifie
\`a l'espace des configurations planaires
d'un $(m+1)$-gone de longueurs $(a_1,\dots ,a_m,|b|)$ modulo les isom\'etries
directes de $\bbr^2$. La classification de ces espaces
a r\'ecemment \'et\'e pouss\'ee jusqu'\`a  $m\leq 8$
dans \cite{HR} (l'espace $M_b$ correspond \`a
$\nua{m+1}{2}(a_1,\dots ,a_m,|b|)$ pour les notations
de \cite{HR}). En g\'en\'eral, si $b\not= 0$, le groupe
$SO(d-1)$ (stabilisateur de $b$ pour l'action de $SO(d)$ sur $\bbr^d$)
agit sur $M_b$ avec quotient hom\'eomorphe \`a
$\nua{m+1}{d}(a_1,\dots ,a_m,|b|)$.

Pour le bras articul\'e $f^a$,
notons $\Delta^a:=\Delta^{f^a}$
la distribution d\'etermin\'ee sur $(S^{d-1})^m$ par $f^a$.
Comme indiqu\'e ci-dessus, le premier travail sera
d'\'etudier les composantes connexes par $\Delta^a$-courbe.
Notre th\'eor\`eme principal ci-dessous les identifie
avec les orbites d'une action sur $(S^{d-1})^m$ par un groupe de Lie,
action que nous allons d\'ecrire maintenant.

Soit $\sharp (a)$ le nombre de valeurs non-nulles de la suite $a_i$~:
$\sharp (a):=\sharp\{a_i\mid i=1\dots ,m,\ a_i\not=0\}.$
Par exemple, $\sharp (3,3,0,1)=2$.
D\'esignons par $\moeb{d-1}$ le groupe des transformations
conformes de $S^{d-1}$ (groupe des transformations de M\"obius,
engendr\'e par les compositions
d'un nombre pair d'inversions, dans $\bbr^d\cup\{\infty\}$,
par rapport \`a des sph\`eres perpendiculaires \`a $S^{d-1}$).
Soient $b_1,\dots ,b_{\sharp(a)}$ les valeurs de la suite $a_i$.
Pour $j\in\{1,\dots ,\sharp(a)\}$, soit
$\moeb{d-1}^{(j)}$  une copie du groupe
$\moeb{d-1}$, que l'on consid\`ere comme agissant
sur $(S^{d-1})^m$ de la mani\`ere
suivante : si $R\in \moeb{d-1}^{(j)}$, on d\'efinit
$R\cdot (z_1,\dots ,z_m):=(z_1',\dots ,z'_m)$ o\`u
$$z'_i:=\left\{\begin{array}{lllll}
R\cdot z_i & \hbox{si $a_i=b_j$}\\
z_i & \hbox{sinon}.
\end{array}\right.$$
Si $j\not=j'$, les actions de $\moeb{d-1}^{(j)}$ et de $\moeb{d-1}^{(j')}$
commuttent et d\'efinissent ainsi une action diff\'erentiable de
$G(a,d):=\moeb{d-1}^{\sharp(a)}$ sur $(S^{d-1})^m$, qui d\'epend de $a$.
Notre th\'eor\`eme principal est~:
\sk{2}\noindent\bf Th\'eor\`eme principal \ \sl
Deux points $z$ et $z'$ de $(S^{d-1})^m$ sont reli\'es par une $\Delta^a$-courbe
si et seulement si ils sont dans la m\^eme orbite pour l'action de
$G(a,d)$. \sk{2}\rm

Le th\'eor\`eme principal sugg\`ere que l'on peut trouver les rel\`evements
horizontaux par des \'equations diff\'erentielles dans le groupe $G(a,d)$.
Ceci sera l'objet d'un travail ult\'erieur.

La preuve du th\'eor\`eme principal occupe les paragraphes
\ref{brasplan} \`a \ref{PrTP}. Elle utilise le th\'eor\`eme de
Sussmann cit\'e auparavant ainsi que le th\'eor\`eme de Chow
(\cite[\S\,0.4]{Gr}, \cite[th. 2.4 et \S\,2.5]{Be},
\cite[Ch. 2]{Mo}). Dans le
paragraphe \ref{Pconseq}, on tire quelques cons\'equences
du th\'eor\`eme principal, par exemple, le fait que les
invariants classiques des transformations de M\"obius,
comme le birapport de $4$ points, produisent des
quantit\'es qui sont conserv\'ees le long des $\Delta^a$-courbes.

Les paragraphes \ref{Phologen} et \ref{S:orbholo} sont consacr\'es
\`a l'\'etude des groupes d'holonomie de $\Delta^a$.
Si $U$ est un ouvert du plan form\'e de valeurs
r\'eguli\`eres du bras  $f^a:(S^{d-1})^m\to \bbr^d$ et si $b\in U$, le groupe
d'holonomie $\calh_b(U)$ est form\'e des diff\'eomorphismes de
$M_b:=f^\mun(b)$ obtenus par transport horizontal au dessus
courbes dans $U$ lisses par morceau qui sont des
lacets en $b$. On montre que pour un
bras g\'en\'erique ($\sharp(a)=m$), le groupe $\calh_b(U)$
agit transitivement sur chaque composante connexe de $M_b$.
Ceci contraste avec le cas $a=(1,\dots,1)$ o\`u l'on montre que,
pour les bras planaires, les orbites de l'action du groupe
d'holonomie sont les lignes de gradient d'une certaine fonction de
Morse $\rho_b:M_b\to S^1$ dont on d\'etermine les points critiques
et leur indice.

Cette recherche a b\'en\'efici\'e du soutien du Fonds National Suisse de la
Recherche Scientifique. L'auteur tient \`a remercier P. de la Harpe
et D. Thurston pour d'utiles conversations.

\section{Bras articul\'es planaires}\label{brasplan}

Ce paragraphe pr\'esente quelques outils pour
la d\'emonstration du th\'eor\`eme principal dans le cas $d=2$
(bras planaires). On retrouvera essentiellement les
m\^emes \'enonc\'es pour les bras dans $\bbr^d$ (\S\,\ref{brasrd}),
dont les d\'emonstrations utiliseront le cas $d=2$ tra\^\i t\'e ici.
De plus, le cas $d=2$ a l'avantage d'utiliser
des techniques plus \'el\'ementaires et des calculs explicites.

Soit $a:=(a_1,\dots ,a_m)\in\bbr^m_{\geq 0}$.
Le bras articul\'e planaire $f^a$ peut s'exprimer comme l'application
$f^a: T^m\to\bbc$ du tore plat $T^m=(S^1)^m$ \`a valeurs complexes~:
$$f^a(z_1,\dots ,z_m)=\sum_{i=1}^m a_iz_i.$$
En param\'etrant le cercle $S^1$ par $q\mapsto e^{iq}$ ($q\in\bbr$),
et en identifiant $\bbc$ avec $\bbr^2$, on regardera aussi $f^a$ comme 
l'application de $\bbr^m\to\bbr^2$ donn\'ee, pour 
$q:=(q_1,\dots ,q_m)\in \bbr^m$, par
\begin{equation}\label{xaya}
f^a(q):=\big(\sum_{i=1}^m a_i\cos q_i\, ,\sum_{i=1}^m a_i\sin q_i\big)
=(X^a(q),Y^a(q)).
\end{equation}
Par abus de langage, les points de $\bbr^m$ seront aussi
appel\'es des {\it configurations}.
La distribution $\Delta^a:=\Delta^{f^a}$ est donc engendr\'ee par les
champs $\grad\, X^a$ et $\grad\, Y^a$. On la regarde, soit comme une distribution
sur $T^m$, soit comme une distribution p\'eriodique sur $\bbr^m$.

Rappelons que le groupe $SU(1,1)$ est form\'e des
($2\times 2)$-matrices complexes $M$ de d\'eterminant 1
telles que $\overline{M}^TJM=J$ o\`u %$J=\pmatrix{1&0\cr 0&-1}$.
$J=\big(\,{}^1_{0}\,{}^{\kern 2pt 0}_{-1}\,\big)$.
Ce sont exactement les matrices de la forme
$\scriptstyle\pmatrix{a& b\cr \bar b& \bar a\cr}$ avec  $\|a\|^2-\|b\|^2=1$.
Le groupe $SU(1,1)$ est le sous-groupe de
$SL(2,\bbc )$ dont l'action sur $\bbc$ par transformations
homographiques pr\'eserve l'orientation et le cercle unit\'e;
il s'identifie donc \`a $\moeb{1}$. Son alg\`ebre
de Lie $su(1,1)$ est form\'ee des matrices $M$ telles que
$J\bar M J = - M^T$ qui s'\'ecrivent donc
$\pmatrix{iu&b\cr \bar b & -iu\cr}$ avec $u\in\bbr$.
Enfin, rappelons que les alg\`ebres de Lie $sl(2,\bbr)$ et $su(1,1)$
sont isomorphes, en fait conjugu\'ees dans $GL(2,\bbc)$~:
$$ Psl(2,\bbr)P^\mun=su(1,1) \ , \hbox{ o\`u }\
P:=\pmatrix{1& -i\cr 1& i\cr}.$$

D\'esignons par $\vect (M)$ l'alg\`ebre de Lie des champs de vecteurs sur une
vari\'et\'e diff\'erentiable $M$. Soit $\call^a$ la sous-alg\`ebre de Lie
de $\vect (T^m)$ (ou de $\vect (\bbr^m)$) engendr\'ee par les
champs $\grad\, X^a$ et $\grad\, Y^a$, o\`u $X^a$ et $Y^a$ sont d\'efinis dans
l'\'equation \eqref{xaya}. 

Soit $\alpha_a : su(1,1)^{\sharp(a)}\to\vect T^m$ 
l'action infinit\'esimale associ\'ee \`a l'action de 
$G(a,2)=SU(1,1)^{\sharp (a)}$ sur $T^m$ d\'efinie dans l'introduction.
Rappelons que $\alpha_a$ est un anti-homomorphisme d'alg\`ebres de Lie.

\begin{Proposition}\label{propab}
L'anti-homomorphisme 
$\alpha_a : su(1,1)^{\sharp(a)}\to\vect (S^1)^m$
est injectif et son image est $\call^a$.
\end{Proposition}

Le reste de ce paragraphe est consacr\'e \`a la preuve de la
proposition \ref{propab}. Nous avons tout d'abord besoin
d'un certain mat\'eriel pr\'eliminaire.
Prenons les notations de l'\'equation \eqref{xaya} et
consid\'erons les champs $S,C\in\vect (\bbr^m)$ donn\'es par
\begin{equation}\label{defcsu}
S:= \pmatrix{\sin q_1\cr\vdots\cr \sin q_m} \qquad  
C:= \pmatrix{\cos q_1\cr\vdots\cr \cos q_m}.
\end{equation} 

On v\'erifie par calcul direct que
\begin{equation} 
[C,S] = \pmatrix{1\cr\vdots\cr 1\cr}=: U \quad , \quad 
[U,C]= -S \quad , \quad  [S,U]= -C.
\end{equation} 
Ces relations sont satisfaites par les g\'en\'erateurs suivants de $su(1,1)$~: 
\begin{equation}\label{tildeSCU}
\tilde C := \pmatrix{0 & i/2 \cr -i/2 & 0\cr} \ , \ 
\tilde S := \pmatrix{0 & 1/2 \cr 1/2 & 0\cr} \ , \ 
\tilde U := \pmatrix{i/2 & 0 \cr 0 & -i/2\cr}.  
\end{equation} 
Les sous-groupes \`a un param\`etre de $SU(1,1)$
correspondants aux deux premiers \'el\'ements sont 
$$Exp(t\, \tilde S) = \pmatrix{\cosh(t/2) & \sinh (t/2)\cr 
\sinh (t/2) & \cosh (t/2)\cr} $$ 
et 
$$Exp(t\, \tilde C) = \pmatrix{\cosh(t/2) & i \sinh (t/2)\cr 
-i \sinh (t/2) & \cosh (t/2)\cr} .$$
%et $$Exp(t\, \tilde U) = \pmatrix{e^{it/2} &  0 \cr 0 & e^{-it/2} \cr} .$$ 
En faisant agir ces sous-groupes \`a un param\`etre sur $S^1$ 
par transformations 
homographiques, on obtient les flots $\Gamma_t^1,\Gamma_t^i:S^1\to S^1$  
($t\in\bbr$) d\'efinis par 
\begin{equation}\label{flotgamma}
\Gamma_t^1(z) :=\frac{\cosh(t/2)\, z + \sinh (t/2)}{\sinh(t/2)\, 
z + \cosh (t/2)}
\end{equation} 
et 
\begin{equation}\label{flotsigma}
\Gamma_t^i(z) :=\frac{\cosh(t/2)\, z + i\sinh (t/2)}{-i\sinh(t/2)\, 
z + \cosh (t/2)}
\end{equation} 
(l'exposant $1$ ou $i$ indique le point fixe stable du flot).
Consid\'erons les champs de vecteurs $S$, $C$, et $U$
de \eqref{defcsu} dans le cas $m=1$, c'est-\`a-dire sur $S^1$.

\begin{Lemme}\label{1param} 
En tant que flot sur $S^1$, on a que

(a)  $\Gamma_t^1$ est
le flot de $-S$, le flot du gradient de la projection
sur l'axe r\'eel.  

(b)  $\Gamma_t^i$ est 
le flot de $C$, le flot du gradient de la projection
sur l'axe imaginaire. 

(c) $z\mapsto e^{it}z$ 
(rotation \`a vitesse constante) est le flot de $U$. 
\end{Lemme} 
 
\preu
Le calcul direct sur la formule \eqref{flotgamma}
montre que $\Gamma_t^1$, d\'efini sur $\bbc$,
est le flot du champ de vecteurs 
\begin{equation}\label{flotasschpvect} 
\frac{d\Gamma_t(z)}{dt}\bigg|_{t=0}= \frac{1-z^2}{2}. 
\end{equation} 
D'autre part, pour $z=e^{i\theta}$ on a  
\begin{equation}\label{normech} 
\begin{array}{lcl} 
\|e^{2i\theta}-1\|^2 &=& \|\cos {2\theta} + i \sin 2\theta -1\|^2=  
(\cos {2\theta}-1)^2 + \sin^2 2\theta = \\ &=& 
2(1-\cos {2\theta}) = 2(1-\cos^2\theta + \sin^2\theta) = 4 \sin^2\theta. 
\end{array}\end{equation} 
Comme $\Gamma_t^1$ pr\'eserve le cercle unit\'e $S^1$,
l'\'equation \eqref{flotasschpvect} 
d\'etermine un champ de vecteurs tangents \`a $S^1$, dont les z\'eros sont $\pm 1$
et qui est  dirig\'e vers la droite. Par \eqref{normech} on d\'eduit que  
$$\frac{d\Gamma_t^1(z)}{dt}\bigg|_{t=0} =
D{\rm expi}(-\sin \theta \frac{\partial}{\partial\theta}).$$ 
Ceci d\'emontre le point a) du lemme \ref{1param}. Le point b)  
se d\'emontre de la m\^eme mani\`ere et le point c) est \'evident.  
\cqfd

Plus g\'en\'eralement, pour $s\in S^1$,
consid\'erons le flot sur $S^1$
\begin{equation}\label{flotgammagen}
\Gamma_t^s(z) :=\frac{\cosh(t/2)\, z + s\sinh (t/2)}{\bar s\sinh(t/2)\,
z + \cosh (t/2)}.
\end{equation}
Soit $g_s:S^1\to\bbr$ la fonction ``$s$-composante''~:
$g_s(z)={\rm Re\,}(z\bar s)$. Soit $\Phi_t^s$
le flot sur $S^1$ de $\grad\, g_s$.
Les points (a) ou (b)
du lemme \ref{1param} sont des cas particuliers
du lemme suivant.

\begin{Lemme}\label{1paramgen}
$\Phi_t^s=\Gamma_t^s$ (\'egalit\'e de flots sur $S^1$).
%Sur $S^1$, le flot $\Gamma_t^s est le flot du gradient de $g_s$.
\end{Lemme}

\preu Soit $m_s:S^1\to S^1$ la multiplication par $s$ et
$m_{\bar s}$ celle par $\bar s$. On a $g_s=g_1\pcirc m_{\bar s}$
et les champs $\grad\,g_1$
et $\grad\,g_s$ sont $m_{\bar s}$-reli\'es~:
$$T_zm_{\bar s}(\grad_z g_s) = \grad_{z\bar s}g_1.$$
On a donc $\Phi_t^s=m_s\pcirc\Phi_t^1\pcirc m_{\bar s}$.
Par le point (a) du lemme \ref{1param}, on a $\Phi_t^1=\Gamma_t^1$,
d'o\`u
$$\Phi_t^s=m_s\pcirc\Gamma_t^1\pcirc m_{\bar s}=\Gamma_t^s. \cqfd$$

Pour $s\in S^1$, soit
$$C_s:=\pmatrix{0 & s/2\cr \bar s/2 & 0\cr}\, \in \, su(1,1).$$
Par exemple, $C_1=\tilde S$ et $C_i=\tilde C$. Comme
$s=e^{i\theta}$, on a que $C_s=P_\theta \tilde S P_{\theta}^\mun$,
o\`u $=P_\theta\in SU(1,1)$ est la matrice diagonale de coefficients
$e^{i\theta/2},e^{-i\theta/2}$.  L'action
par transformation homographique de $P$ sur $S^1$
est la rotation $m_s$.
On a donc que $\exp (t C_s)$ g\'en\`ere sur $S^1$ le flot $\Gamma^s_t$.

Soit $\alpha:su(1,1)\to \vect(S^1)$
l'action infinit\'esimale de l'action de $SU(1,1)$
sur $S^1$. Vu ce qui pr\'ec\`ede,
le lemme \ref{1param} se paraphrase en

\begin{Lemme}\label{1paraminf}
$\alpha(C_s)=\grad\,g_s$.
\cqfd
\end{Lemme}

\begin{Remarque}\rm Comme
$\alpha$ est un anti-homomorphisme d'alg\`ebres de Lie,
on a bien que
$$\alpha(\tilde U)=\alpha([\tilde C,\tilde S])=
-[\alpha(\tilde C),\alpha(\tilde S)])=-[C,-S]=U.$$
\end{Remarque}

\paragraph{Preuve de la proposition \ref{propab}.}
Elle se fait par r\'ecurrence sur $\sharp(a)$.
Le cas $\sharp(a)=0$ est banal car 
$su(1,1)^0=0$ (par convention) et $\call^a=0$ puisque
l'application $f^a$ est constante. L'injectivit\'e
de $\alpha_a$ est aussi \'evidente puisque l'action de 
$G(a,d)$ sur $S^{d-1}$ est effective.

Supposons que $\sharp(a)=1$. Soit $b$ l'unique
valeur positive de la suite $a_i$.
Les composantes non-nulles de $\grad\, X^a$
sont \'egales \`a celles de $-bS$ et celles de
$\grad\, Y^a$ sont \'egales \`a celles de $bC$. De plus,
$\grad\, X^a$ et $\grad\, Y^a$ ont les m\^emes composantes
nulles. Il s'en suit que
l'alg\`ebre de Lie $\call^a$ ne d\'epend pas de $b$
 et qu'elle est
isomorphe \`a l'alg\`ebre de Lie engendr\'ee par
par $S$ et $C$.
L'action de $G(a,2)=SU(1,1)$ sur $T^m$
ne d\'ependant pas non plus de $b$, on peut supposer 
que $b=1$. Dans ce cas, 
pour $X\in su(1,1)$
les composantes non-nulles de $\alpha_a(X)$ 
sont \'egales \`a $\alpha(X)$. Par les parties
(a) et (b) du lemme \ref{1paraminf},
on en d\'eduit que 
$\alpha_a(\tilde S)=\grad\, X^a$
et $\alpha_a(\tilde C)=\grad\, Y^a$.
L'image de $\alpha_a$ contient donc $\call^a$. 
Comme $\tilde S$ et $\tilde C$ 
engendrent $su(1,1)$ en tant qu'alg\`ebre de Lie,
l'image de $\alpha_a$ 
est aussi contenue dans $\call^a$, ce qui d\'emontre 
la proposition \ref{propab} lorsque $\sharp(a)=1$.

Supposons, par hypoth\`ese de r\'ecurrence, que
la proposition \ref{propab} soit vraie pour
$\sharp(a)\leq j-1$. Soit $a$ tel que $\sharp(a)=j$.
Pour $k\in\bbn$, introduisons les champs de vecteurs
suivants sur $\bbr^m$ 
\begin{equation}\label{chpsak}
A^k:=\pmatrix{a_1^k\cr\vdots\cr a_m^k\cr} \ ,\ 
A^kS:=\pmatrix{a_1^k\sin q_1\cr\vdots\cr a_m^k\sin q_m\cr}\ ,\ 
A^kC:=\pmatrix{a_1^k\cos q_1\cr\vdots\cr a_m^k\cos q_m\cr}.
\end{equation}
L'espace vectoriel engendr\'e par les champs
\begin{equation}\label{Ak}
A^{2k-1}C\ ,\ A^{2k-1}S\ \hbox{ et }\ A^{2k}\qquad (k\geq 1)
\end{equation}
est une sous-alg\`ebre de Lie de $\vect(\bbr^m)$. En effet, on a
$$[A^{2k-1}C,A^{2l-1}S]=A^{2(k+l-1)}\quad ,\quad
[A^{2k},A^{2l-1}C]=-A^{2(k+l)-1}S$$ $$
[A^{2k-1}S,A^{2l}]=-A^{2(k+l)-1}C  \quad ,\quad [A^{2k},A^{2l}]=0.$$
Cette alg\`ebre de Lie
est engendr\'ee par $\grad\, X^a = -AS$ et $\grad\, Y^a = AC$,
donc c'est l'alg\`ebre de Lie $\call^a$.

Supposons que $a_m\not= 0$ (le
raisonnement est analogue si c'est un autre $a_i$). Consid\'erons 
les projecteurs lin\'eaires $\pi',\pi'' : \bbr^m\to\bbr^m$ d\'efinis par
$\pi'(x_1,\dots ,x_m):=(x_1',\dots ,x_m')$ et
$\pi''(x_1,\dots ,x_m):=(x_1'',\dots ,x_m'')$ o\`u
$$x'_i:=\left\{\begin{array}{lll}
x_i & \hbox{si $a_i=a_m$}\\
0   &  \hbox{sinon}
\end{array}\right. \ ,\ \
x''_i:=\left\{\begin{array}{lll}
0 & \hbox{si $a_i=a_m$}\\
x_i   &  \hbox{sinon}
\end{array}\right.$$
On consid\'erera aussi $\pi'$ et $\pi''$ comme des morphismes
du fibr\'e tangent $T\bbr^m$ sur lui m\^eme, au dessus de l'identit\'e
de $\bbr^m$, en appliquant $\pi'$ ou $\pi''$ sur chaque fibre.
Notons $a':=\pi'(a)$ et $a'':=\pi''(a)$. Par exemple, si $m=5$ et
$a=(2,1,3,1,1)$, alors $a'=(0,1,0,1,1)$ et $a''=(2,0,3,0,0)$.

On a $\pi'(\grad\, X^a)=\grad\, X^{a'}$, $\pi'(\grad\, Y^a)=\grad\, Y^{a'}$,
etc, ce qui montre que $\Delta^a=\Delta^{a'}\oplus\Delta^{a''}$.
De plus, la forme particuli\`ere des champs de \eqref{chpsak} 
et des actions fait que $\call^a=\call^{a'}\oplus \call^{a''}$.
Bref, tout se passe composante par composante et il en est
de m\^eme des actions~: $\alpha:\lie(G(a,2))\to\vect(T^m)$
se d\'ecompose en $\alpha_a=\alpha_{a'}\oplus\alpha_{a''}$.
Comme $\sharp(a')=1$ et $\sharp(a'')=\sharp(a)-1$, on aura,
par hypoth\`ese de r\'ecurrence, que l'image de $\alpha_{a'}$
est $\call_{a'}$ et que l'image de $\alpha_{a''}$
est $\call_{a''}$. On en d\'eduit que l'image de $\alpha_a$
est bien $\call^a$, ce qui ach\`eve la d\'emonstration 
de la proposition \ref{propab}. \cqfd

\begin{Corollaire}\label{deuxgen}
En tant qu'alg\`ebre de Lie, $su(1,1)^m$ est,
quelque soit $m$, engendr\'ee par 2 \'el\'ements.
\end{Corollaire}

\preu 
Soit $a=(a_1,\dots,a_m)\in\bbr^m$ avec $0<a_1<\cdots <a_m$
(donc $\sharp(a)=m$).
Par la proposition \ref{propab},
$\alpha_a:su(1,1)^m\to\call^a$ est un
(anti)-isomorphisme d'alg\`ebres de Lie. Or, $\call^a$
est engendr\'ee, en tant qu'alg\`ebre de Lie, par 
$\grad\, X^a$ et $\grad\, X^a$.\cqfd

\begin{Remarque}\label{deuxgenrem} \rm
Pour $m=\sharp(a)$, les deux g\'en\'erateurs de $su(1,1)^m$ du corollaire
\ref{deuxgen} sont
$$A\tilde C:=(a_1\tilde C,\dots ,a_m\tilde C)\quad\hbox{ et }\quad
A\tilde S:=(a_1\tilde S,\dots ,a_m\tilde S),$$
avec les $\tilde C,\tilde S\in su(1,1)$ d\'efinis en \eqref{tildeSCU}.
\end{Remarque}

\begin{ex}\rm Les champs $A^1,A^2\dots A^k$ de \eqref{chpsak}
sont lin\'eairement d\'ependants d\`es que $k>\sharp(a)$. 
La relation de d\'ependance lin\'eaire se 
r\'epercute sur les $A^iC$ et $A^iS$, donnant une pr\'esentation
int\'eressante de l'alg\`ebre de Lie $su(1,1)^{\sharp(a)}$.
Par exemple, si $a=(a_1,a_2)$ avec $a_1\not= a_2$,
on obtient une pr\'esentation de $su(1,1)^2$ engendr\'ee comme espace
vectoriel par
$$ AC\, ,\,AS\, ,\,A^2\, ,\,A^3C\, ,\,A^3S\, ,\,A^4$$ 
et les crochets non-nuls sont :
$$\begin{array}{cclccclcccl} 
[AC,AS] &=& A^2  \\[2pt]  
[A^2,AC] &=& -A^3S \\[2pt]  
[A^2,AS] &=& A^3C \\[2pt] 
[AC,A^3S]  &=& A^4 \\[2pt]  
[A^2,A^3C] &=& [A^4,AC] = a_1^2a_2^2AS-(a_1^2+a_2^2)A^3S\\[2pt]  
 [A^2,A^3S] &=& [A^4,AS] = -a_1^2a_2^2AC+(a_1^2+a_2^2)A^3C \\[2pt] 
[A^3C,A^3S] &=& -a_1^2a_2^2A^2+(a_1^2+a_2^2)A^4  \\[2pt]  
[A^4,A^3C] &=& a_1^2a_2^2(a_1^2+a_2^2)AS - (a_1^4+a_1^2a_2^2+a_2^4)A^3S 
 \\[2pt]  
[A^4,A^3S] &=& -a_1^2a_2^2(a_1^2+a_2^2)AC +  
(a_1^4+a_1^2a_2^2+a_2^4)A^3C .
\end{array}$$
Les quatre premi\`eres relations montrent que $su(1,1)^2$ est,
en tant qu'alg\`ebre de Lie, engendr\'ee par $AC$ et $AS$.
\end{ex}

\section{Bras articul\'es dans $\bbr^d$}\label{brasrd}

Dans ce paragraphe, nous g\'en\'eralisons aux bras dans 
$\bbr^d$ les r\'esultats obtenus,
dans le paragraphe \ref{brasplan}, pour les bras planaires.
Les d\'emonstrations utilisent les r\'esultats du \S\,\ref{brasplan}.

Soit $s\in S^{d-1}$. Consid\'erons le flot
$\Gamma_t^s$ sur $S^{d-1}$ form\'e des
transformations de M\"obius d\'etermin\'ees par
la condition suivante~:
si $\varphi$ est une transformation de M\"obius de
$\bbr^d\cup\{\infty\}$ telle que $\varphi(-s=0$ et
$\varphi(s=\infty$, alors
$\varphi\pcirc\Gamma_t^s\pcirc\varphi^\mun(x)=e^t\,x$, pour
tout $x\in\varphi(S^{d-1})$.
En terme d'isom\'etries hyperboliques du
disque de Poincar\'e $D^d$,
$\Gamma_t^s$ correspond \`a un flot de translations
hyperboliques laissant invariant chaque plan qui contient
la g\'eod\'esique allant de $-s$ \`a $s$.

Soit $g_s:S^{d-1}\to\bbr$ la fonction
$g_s(x)=\llangle{x}{s}$, o\`u $\llangle{}{}$ d\'enote
le produit scalaire standard dans $\bbr^d$.
Soit $\Phi_t^s$ le flot, sur $S^{d-1}$, du gradient
de $g_s$.

\begin{Lemme}\label{1paramd}
$\Gamma_t^s=\Phi_t^s$ (\'egalit\'e de flots sur $S^{d-1}$).
\end{Lemme}

\preu
Les deux flots ayant $\pm s$ comme points fixes,
il suffit de v\'erifier que $\Gamma_t^s(y)=\Phi_t^s(y)$
pour $y\in S^{d-1}$ diff\'erent de $\pm s$.

Soit $\Pi$ le sous-espace vectoriel de $\bbr^d$
engendr\'e par $y$ et $s$. Soit $h:\bbr^d\to\bbr^d$
la r\'eflexion orthogonale par rapport au plan $\Pi$.
La relation g\'en\'erale entre $\grad(g_s\pcirc h)$
et $\grad\, g_s$ est
\begin{equation}\label{hgrad}
(D_xh)^*(\grad_{h(x)}g_s) = \grad_{x}(g_s\pcirc h),
\end{equation}
o\`u $(D_xh)^*$ est l'adjoint de $D_xh$. Comme $h$ est
une isom\'etrie lin\'eaire, on a $(D_xh)^*=h^*=h^\mun=h$.
D'autre part, $g_s\pcirc h=g_s$. Enfin, pour $x\in\Pi$,
on a $x=h(x)$ et la formule \eqref{hgrad} se r\'eduit \`a
$h(\grad_{x}g_s)=\grad_{x}g_s$. Cela prouve que,
si $x\in\Pi$, alors $\grad_{x}g_s\in\Pi$, ce qui
implique que $\grad_x g_s=\grad_x (g_s)_{|\Pi}$.

Le flot $\Phi_t^s$ pr\'eserve donc le plan $\Pi$. Comme
c'est aussi le cas du flot $\Gamma_t^s$,
il suffit de montrer leur \'egalit\'e sur le cercle
$\Pi\cap S^{d-1}$. Soit $\mu:\bbc\to\Pi$ une isom\'etrie
euclidienne $\bbr$-lin\'eaire telle que $\mu(1)=s$
(il en existe deux, mais ce choix est sans importance).
On a, sur $\Pi\cap S^{d-1}$, que
$\mu\pcirc\Gamma_t^s\pcirc\mu^\mun=\Gamma_t^1$, o\`u
$\Gamma_t^1$ est le flot d\'efini dans l'\'equation
\eqref{flotgamma}. De m\^eme,
$\mu\pcirc\Phi_t^s\pcirc\mu^\mun$ est le flot
sur $S^1$ du gradient de la projection sur
l'axe r\'eel. Le fait que $\Gamma_t^s=\Phi_t^s$
sur $\Pi\cap S^{d-1}$ est exactement le contenu de la partie
(a) du lemme \ref{1param}. \cqfd

On peut \'evidemment voir $\Gamma_t^s$
comme un sous-groupe \`a un param\`etre de $\moeb{d-1}$.
Soit $C_s$ l'\'el\'ement de $\lie(\moeb{d-1})$
tel que $\exp(tC_s)=\Gamma_t^s$.
Soit $\alpha:\lie(\moeb{d-1})\to \vect(S^{d-1})$
l'action infinit\'esimale de l'action de $\moeb{d-1}$
sur $S^{d-1}$. Le lemme \ref{1paramd} se paraphrase en

\begin{Lemme}\label{1paraminfd}
$\alpha(C_s)=\grad\, g_s$. \cqfd
\end{Lemme}

Soit $s'\in S^{d-1}$ avec $\llangle{s}{s'}=0$.
D\'esignons par $\Pi(s,s')$ le 2-plan de $\bbr^d$
engendr\'e par $s$ et $s'$ et orient\'e par cette base.
On appelle {\it rotation d'angle $t$ dans $\Pi(s,s')$}
l'\'el\'ement $R\in SO(d)$ tel que
$$R(s)=\cos t s + \sin t s' \ , \
R(s')=-\sin t s + \cos t s'$$
et $R(x)=x$ pour tout $x$ dans le compl\'ement
orthogonal de $\Pi(s,s')$.

\begin{Lemme}\label{crochets}
Le sous-groupe \`a un param\`etre $\exp(t[C_s,C_{s'}])$
est la rotation d'angle $-t$ dans le plan $\Pi(s,s')$.
\end{Lemme}

\preu
Soit $\Pi=\Pi(s,s')$.
Comme $\Gamma_t^s$ et $\Gamma_t^{s'}$ pr\'eservent
le plan $\Pi$, les champs de vecteurs
$\alpha(C_s)$ et $\alpha(C_{s'})$ sont tangents \`a $\Pi$.
Le crochet
$[\alpha(C_s),\alpha(C_{s'})]=-\alpha([C_s,C_{s'}])$ est
donc tangent \`a $\Pi$ et le sous-groupe \`a un param\`etre
$\exp(t[C_s,C_{s'}])$ pr\'eserve le plan $\Pi$.

Soit $\mu:\bbc\to\Pi(s,s')$ l'isom\'etrie
euclidienne $\bbr$-lin\'eaire telle que $\mu(1)=s$
et $\mu(i)=s'$. On a vu dans la preuve du lemme
\ref{1paramd} que $\mu$ conjugue $C_s$ et
$C_{s'}$ avec, respectivement, les \'el\'ements
$\tilde S$ et $\tilde C$ de $su(1,1)$
d\'efinis par les \'equations \eqref{tildeSCU}.
L'isom\'etrie $\mu$ conjugue donc $[C_s,C_{s'}]$
avec $[\tilde S,\tilde C]=-\tilde U$
dont le flot sur $S^1$ est
la rotation d'angle $-t$ (partie (c)
du lemme \ref{1param}). Cela montre
que le sous-groupe \`a un param\`etre $\exp(t[C_s,C_{s'}])$
se restreint \`a $\Pi$ en la rotation d'angle $-t$.

Il reste \`a montrer que si $x\in S^{d-1}$ est
orthogonal \`a $s$ et \`a $s'$, alors
$\exp(t[C_s,C_{s'}])\,x=x$. Soit $S$ la $2$-sph\`ere
de rayon 1 dans l'espace vectoriel engendr\'e par
$s$, $s'$ et $x$. Comme $\Gamma_t^s=\exp(tC_s)$
et $\Gamma_t^{s'}=\exp(tC_{s'})$
pr\'eservent $S$, ainsi en est-il de
$\exp(t[C_s,C_{s'}])$.
Comme le flot $\exp(t[C_s,C_{s'}])$ pr\'eserve $\Pi$, il
pr\'eserve la famille des (grands) cercles
de $S$ qui sont perpendiculaires \`a $\Pi$.
Leur intersection commune $\{\pm x\}$ est donc
pr\'eserv\'ee par $\exp(t[C_s,C_{s'}])$. Par continuit\'e,
on en d\'eduit que $\exp(t[C_s,C_{s'}])\,x=x$. \cqfd

\begin{Lemme}\label{engd}
Les \'el\'ements $C_s$, pour $s\in S^{d-1}$,
engendrent $\lie(\moeb{d-1})$ en tant qu'alg\`ebre de Lie.
\end{Lemme}

\preu
Soit $L$ la sous-alg\`ebre de Lie de $\lie(\moeb{d-1})$
engendr\'ee par les \'el\'ements $C_s$. Il suffit de voir que
le groupe $G=\exp(L)$ est \'egal \`a $\moeb{d-1}$.

Identifions $\moeb{d-1}$ avec ${\rm Iso}^+(D^d)$, le groupe
des isom\'etries hyperboliques du disque de Poincar\'e $D^d$
qui pr\'eservent l'orientation. L'orbite de l'origine
par $\exp(tC_s)$ est le segment allant de $-s$ \`a $s$.
On en d\'eduit que $G$ agit transitivement sur $D^d$.
Il suffit donc de montrer que $G$ contient le stabilisateur
de l'origine qui est $SO(d)$. Mais, par le lemme \ref{crochets},
$G$ contient toutes les rotations dans tous les $2$-plans
de $D^d$ et ces rotations engendrent $SO(d)$. \cqfd

Nous allons maintenant d\'emontrer la g\'en\'eralisation
de la proposition \ref{propab}.
Soit $a=(a_1,\dots,a_m)\in (\bbr_{\geq 0})^m$.
Si $f^a$ est un bras articul\'e de type $a$,
on d\'esigne par $g_s^a:(S^{d-1})^m\to\bbr$ la fonction
d\'efinie par $g_s^a(z)=\llangle{f^a(z)}{s}$.
Soit $\call^a$ la sous-alg\`ebre de Lie
de $\vect (S^{d-1})^m$ engendr\'ee par les
champs $\grad\, g_s^a$ pour tout $s\in S^{d-1}$.
Soit $\alpha_a : \lie(G(a,d))\to\vect (S^{d-1})^m$
l'anti-homomorphisme d'alg\`ebre de Lie associ\'e
\`a l'action de $G(a,d)$ d\'efinie dans l'introduction.

\begin{Proposition}\label{propabd}
L'anti-homomorphisme
$\alpha_a : \lie(G(a,d))\to\vect (S^{d-1})^m$
est injectif et son image est $\call^a$.
\end{Proposition}

\preu Comme pour la proposition \ref{propab}, la
d\'emonstration se fait par r\'ecurrence sur $\sharp(a)$.
Le cas $\sharp(a)=0$ est banal. L'injectivit\'e
de $\alpha_a$ est aussi \'evidente puisque l'action de
$G(a,d)$ sur $S^{d-1}$ est effective.

Supposons que $\sharp(a)=1$. Soit $b$ l'unique
valeur positive de la suite $a_i$ et soit
$\bar a = \frac{1}{b}\, a$. Comme
$g_s^a(z)=\sum_{j=1}^m a_j\llangle{z_j}{s}$,
%la $j$-\`eme
%composante de $\grad\,g_s^a$ dans $(\bbr^d)^m$ est
%\'egale \`a $a_j\,\grad\, g_s$. On a donc que
on a
$\grad\,g_s^a= b\,\grad\,g_s^{\bar a}$, d'o\`u il
r\'esulte que $\call^a=\call^{\bar a}$. On a aussi
que $G(\bar a,d)=G(a,d)$ et $\alpha_{\bar a}=\alpha_a$.

On peut donc supposer que $b=1$. Dans ce cas,
pour $X\in\lie(G(a,d))=\lie(\moeb{d-1})$,
les composantes non-nulles de $\alpha_a(X)$
sont \'egales \`a $\alpha(X)$. Par le lemme \ref{1paraminfd},
$\alpha(C_s)=\grad\, g_s$, donc
$\alpha_a(C_s)=\grad\, g_s^a$ et
l'image de $\alpha_a$ contient $\call^a$.
D'autre part,
par le lemme \ref{engd}, les \'el\'ements $C_s$
engendrent $\lie(\moeb{d-1})$ comme alg\`ebre de Lie,
ce qui entra\^\i ne que l'image de $\alpha_a$
est contenue dans $\call^a$. Le cas $\sharp(a)=1$
est ainsi d\'emontr\'e.

Le pas de r\'ecurrence est exactement le m\^eme que
pour la proposition \ref{propab}. \cqfd

\begin{Corollaire}\label{deuxgend}
En tant qu'alg\`ebre de Lie, $\lie(\moeb{d-1})^m$ est,
quelque soit $m$, engendr\'ee par $d$ \'el\'ements.
\end{Corollaire}

\preu
Soit $a=(a_1,\dots,a_m)\in\bbr^m$ avec $0<a_1<\cdots <a_m$
(donc $\sharp(a)=m$).
Par la proposition \ref{propabd},
$\alpha_a:\lie(\moeb{d-1})^m\to\call^a$ est un
(anti)-isomorphisme d'alg\`ebres de Lie. Or, $\call^a$
est engendr\'ee, en tant qu'alg\`ebre de Lie, par
les $d$ champs $\grad\, g_s^a$ pour $s\in S^{d-1}$ parcourant
les \'el\'ements d'une base de $\bbr^d$. \cqfd

\section{Preuve du th\'eor\`eme principal}\label{PrTP}

Nous allons tout d'abord rappeler un th\'eor\`eme de
Sussmann \cite[Th.\,4.1]{Su} sous la forme
simplifi\'ee dont nous aurons besoin ici.
Soit $M$ une vari\'et\'e diff\'erentiable et $\cald\subset\vect (M)$.
Soit $\Delta_\cald$ la distribution sur $M$ engendr\'ee par
$\cald$, c'est-\`a-dire que $(\Delta_\cald)_x$
est le sous-espace vectoriel de
$T_xM$ engendr\'e par $A_x$ pour $A\in\cald$.
Supposons que les champs $A\in\cald$ sont complets,
c'est-\`a-dire dont leur flot $\Phi^A_t$ est d\'efini
pour tout $t\in\bbr$. Soit $H$
le sous-groupe de diff\'eomorphismes de $M$
engendr\'e par les $\Phi^A_t$ pour tout $A\in\cald$ et $t\in\bbr$.
Soit $P_\cald$ la plus petite distribution sur $M$
qui contient $\Delta_\cald$ et qui est invariante par l'action de $H$.

Soit $K$ une distribution sur $M$.
Rappelons qu'une {\it vari\'et\'e int\'egrale} de $K$
est une sous-vari\'et\'e $V$ de $M$ telle
que $T_xV=K_x$ pour tout $x\in V$. La distribution
$K$ est {\it int\'egrable} si tout point de $M$ appartient \`a une
vari\'et\'e int\'egrable de $K$. La relation
``\^etre dans une m\^eme vari\'et\'e int\'egrale'' engendre alors une
relation d'\'equivalence sur $M$ dont les classes d'\'equivalence sont
les {\it feuilles} d'un feuilletage. Ces
feuilles ne sont pas, en g\'en\'eral, toutes de m\^eme dimension;
l'exemple standard est donn\'e
par les orbites d'une action d'un groupe de Lie
sur une vari\'et\'e. Le th\'eor\`eme de Sussmann \cite[Th.\,4.1]{Su},
lorsque les champs de $\cald$ sont complets (c'est le cas si,
par exemple, si $M$ est compacte), prend la forme suivante~:

\begin{Theorem}[Sussmann]
Soit $\cald\subset\vect (M)$ form\'e de champs complets.
La distribution $P_\cald$ est int\'egrable et ses feuilles
sont les orbites de l'action de~$H$.
\end{Theorem}

Le corollaire suivant est
bien connu des sp\'ecialistes (voir \cite[p.230]{Mo}) mais ne semble
pas \^etre \'enonc\'e explicitement dans la litt\'erature.

\begin{Corollaire}\label{cosussmann}
Soit $\cald\subset\vect (M)$ form\'e de champs complets.
Soient $x,y\in M$.
Les deux conditions suivantes sont \'equivalentes~:
\begin{enumerate}
\renewcommand{\labelenumi}{(\roman{enumi})}
\item $x$ et $y$ sont joignables par une $\Delta_\cald$-courbe.
\item $x$ et $y$ sont joignables par une courbe qui est
une succession de trajectoires de champs de $\cald$.
\end{enumerate}
\end{Corollaire}

\preu
Il est banal que (ii) entra\^\i ne (i). R\'eciproquement,
soit $c$ une $\Delta_\cald$-courbe joignant $x$ \`a $y$.
Comme la distribution $P_\cald$ contient $\Delta_\cald$,
la courbe $c$ est aussi une $P_\cald$-courbe. Elle est
donc contenue dans une feuille de $P_\cald$. Par le
th\'eor\`eme de Sussmann, les points $x$ et $y$ sont dans
la m\^eme orbite de~$H$, ce qui est \'equivalent \`a
la condition (ii). \cqfd

Nous pouvons maintenant entreprendre la d\'emonstration
du th\'eor\`eme principal.
Soit $G:=G(a,d)$.
Soient $x,y$ deux points de $(S^{d-1})^m$
joignables par une $\Delta^a$-courbe.
Rappelons que $\Delta^a=\Delta_\cald$
pour $\cald=\{\grad\, g^a_s\mid s\in S^{d-1} \}$.
Les champs $\grad\, g^a_s$ sont complets puisque
$(S^{d-1})^m$ est compacte. Par le corollaire
\ref{cosussmann}, il existe une courbe reliant $x$ \`a $y$ qui est
une succession de trajectoires de champs de $\cald$.
Par la proposition \ref{propabd}, les \'el\'ements de $\cald$
sont des champs fondamentaux de
l'action de $G$ sur $(S^{d-1})^m$.
Les points $x$ et $y$ sont dans une m\^eme orbite pour l'action
de $G$.

R\'eciproquement, soit $z\in (S^{d-1})^m$.
Consid\'erons l'application diff\'erentiable
$\beta : G\to (S^{d-1})^m$ donn\'ee par $\beta(g):=g\, z$,
dont l'image est l'orbite de~$z$.

Soient $e_1,\dots,e_d\in S^{d-1}$ une base de $\bbr^d$.
Par la proposition \ref{propabd}, il existent
$d$ \'el\'ements $C^a_{e_i}\in\lie(G(a,d))$ tels que
$\alpha_a(C^a_{e_i})=\grad\,g^a_{e_i}$. Ces champs
$C^a_{e_i}$ engendrent $\lie(G(a,d))$ en tant qu'alg\`ebre de Lie
(voir la preuve du corollaire~\ref{deuxgend}).

La courbe
$t\mapsto {\rm exp\,}(tC^a_{e_i})\, g\, z$ dans $(S^{d-1})^m$
repr\'esente donc le vecteur tangent
$\grad_{gz}g^a_{e_i}$.
Le champ sur $G$ invariant \`a droite $\underline{C^a_{e_i}}$
 correspondant \`a $C^a_{e_i}$
est donc $\beta$-reli\'e au champ
$\grad\, g^a_{e_i}$. Comme les $C^a_{e_i}$
engendrent $\lie(G(a,d))$
en tant qu'alg\`ebre de Lie, les champs $\underline{C^a_{e_i}}$
 engendrent
${\rm Lie}^-(G)$, l'alg\`ebre de Lie des champs invariants \`a droite sur $G$.
Comme les valeurs des champs de ${\rm Lie}^-(G)$ engendrent
l'espace tangent $T_gG$ en tout $g\in G$, le th\'eor\`eme de Chow
(\cite[\S\,0.4]{Gr}, \cite[th. 2.4 et \S\,2.5]{Be}, \cite[Ch. 2]{Mo})
implique que $G$ est connexe par courbes qui sont des
successions de trajectoires des champs $\underline{C^a_{e_i}}$.
L'image de telles courbes par $\beta$
\'etant des $\Delta^a$-courbes, on en d\'eduit que l'orbite $G\,z$
de $z$ est connexe par $\Delta^a$-courbes.

\begin{Remarque} Changement de la métrique riemannienne. \rm
Le fibré tangent à $(S^{d-1})^m$ se décompose en somme directe
$T(S^{d-1})^m=(TS^{d-1})^m$. Si l'on munit $(S^{d-1})^m$
de la métrique riemanienne pour laquelle
$\|(v_1,\dots,v_m)\|=\sum_{j=1}^m a_j\|v_j\|$,
la distribution $\Delta^a$ pour cette nouvelle métrique est
égale à celle, dans la métrique standard, pour le bras de type
$a=(1,\dots,1)$ et $G(a,d)$ est alors isomorphe à $\moeb{d-1}$.
\end{Remarque}

\section{Cons\'equences du th\'eor\`eme principal}\label{Pconseq}

Soit $f^a : (S^{d-1})^m\to\bbr^d$ un bras articul\'e, avec $a=(a_1,\dots,a_m)$.
La proposition \ref{restentalignes} ci-dessous est une cons\'equence
banale du th\'eor\`eme principal.

\begin{Proposition}\label{restentalignes}
Supposons qu'il existe $i,j$ avec $a_i=a_j$. Soit
$z(t)\in (S^{d-1})^m$, $t\in [0,1]$, une $\Delta^a$-courbe.
Si $z_i(0)=z_j(0)$, alors $z_i(t)=z_j(t)$ pour tout $t\in [0,1]$.\cqfd
\end{Proposition}

Le th\'eor\`eme principal entra\^\i ne que les invariants par
transformations de M\"obius, appliqu\'es aux $z_i$ correspondants
\`a des segments de m\^eme longueur, sont invariants par $\Delta^a$-courbes.
Cela donne toute une famille d'int\'egrales premi\`eres
(constantes du mouvement) des $\Delta^a$-courbes.

Par exemple, si $d=2$,
trois points  $z_1,z_2,z_3\in S^1$ qui sont deux \`a deux distincts
d\'eterminent une orientation $\calo(z_1,z_2,z_3)$ de $S^1$, d\'efinie comme
le sens de parcours tel que le chemin allant de $z_1$ \`a $z_3$
passe par $z_2$. Pour
quatre $z_i\in S^1$ distincts, on a le birapport
$$b(z_i,z_j,z_k,z_l)=\frac{z_i-z_k}{z_i-z_l}\,\frac{z_j-z_l}{z_j-z_k}\,\in\bbr$$
(ce birapport est r\'eel puisque les quatre
$z_i$ sont sur un m\^eme cercle).

Pour $d=3$ et quatre $z_i\in S^2$ distincts, on a le {\it birapport
complexe} $b_\bbc(z_i,z_j,z_k,z_l)\in\bbc$ d\'efini de la fa\c con suivante.
On choisit une transformation de M\"obius $\sigma$ de $\bbr^3\cup\{\infty\}$
qui envoie $S^2$ sur le plan horizontal, identifi\'e avec $\bbc$,
telle que $\sigma : S^2\to\bbc$ pr\'eserve l'orientation. On peut
alors d\'efinir
$$b_\bbc(z_i,z_j,z_k,z_l):=
\frac{\sigma(z_i)-\sigma(z_k)}{\sigma(z_i)-\sigma(z_l)}\,
\frac{\sigma(z_j)-\sigma(z_l)}{\sigma(z_j)-\sigma(z_k)}\,\in\bbc.$$

Enfin, pour $d>3$, on a le {\it birapport faible}
$$b_+(z_i,z_j,z_k,z_l):=
\frac{|z_i|-|z_k|}{|z_i|-|z_l|}\,\frac{|z_j|-|z_l|}{|z_j|-|z_k|}\,
\in \bbr_{>0},$$
qui est un invariant de $\moeb{d-1}$ (voir \cite[p. 32]{Bn}).
Ces invariances donnent imm\'ediatement la proposition suivante.

\begin{Proposition}\label{birapp}
Soit $f^a:(S^{d-1})^m\to\bbr^d$ un bras articul\'e.

(a) Supposons que $d=2$. Soient $i,j,k$ avec $a_i=a_j=a_k$.
Soient $z(t)\in T^m$, $t\in [0,1]$ une $\Delta^a$-courbe.
Supposons que
et $z_i(0)\not=z_j(0)\not=z_k(0)\not=z_i(0)$.
Alors, pour tout $t$, on a $z_i(t)\not=z_j(t)\not=z_k(t)\not=z_i(t)$
et l'orientation $\calo(z_i(t),z_j(t),z_k(t))$
est constante.

(b) Soient $i,j,k,l$ avec $a_i=a_j=a_k=a_l$.
Soient $z(t)\in T^m$, $t\in [0,1]$ une $\Delta^a$-courbe.
Si les quatre points initiaux $z_i(0), z_j(0)$, etc, sont tous distincts,
ils le restent tout au long de la $\Delta^a$-courbe et
\begin{enumerate}
\item  si $d=2$, le birapport $b(z_i(t),z_j(t),z_k(t),z_l(t))$ est constant.
\item  si $d=3$, le birapport complexe $b_\bbc(z_i(t),z_j(t),z_k(t),z_l(t))$ est constant.
\item  pour tout $d\geq 2$, le birapport faible $b_+(z_i(t),z_j(t),z_k(t),z_l(t))$ est constant.
\cqfd
\end{enumerate}
\end{Proposition}

\begin{ccote}\rm Dans le cas $a=(1,\dots,1)$ et $d=2$ ou $3$, les birapports
permettent d'obtenir des invariants complets pour la connexit\'e
par $\Delta^a$-courbes. La situation est r\'esum\'ee par le tableau suivant.
Un tiret dans la colonne de droite veut dire une condition vide,
donc toujours vraie. Soit donc $f^a:(S^{d-1})^m\to\bbr^d$
un bras articul\'e de type $a=(1,\dots,1)$.
Soient $z^0,z^1\in (S^{d-1})^m$.

\begin{center}
\begin{tabular}{c|c|c} \sl\small
$m$ & $d$ &
\begin{minipage}[c]{90mm} \begin{center} \small\sl
Condition \'equivalente \`a l'existence d'une\\ $\Delta^a$-courbe joignant
$z^0$ \`a $z^1$
\end{center}\end{minipage}
\\[3mm] \hline \rule{0pt}{3ex}
$2$ & $\geq 2$ & -- \\[2mm] \hline  \rule{0pt}{4ex}
$3$ & $2$ &
\begin{minipage}[c]{80mm} \begin{center} \small\sl
%ordre cyclique sur $S^1$ induit par $(z^0_1,z^0_2,z^0_3)$
%\'egal\\ \`a celui induit par $(z^1_1,z^1_2,z^1_3)$
$\calo(z^0_1,z^0_2,z^0_3)=\calo(z^1_1,z^1_2,z^1_3)$
\end{center}\end{minipage} \\[3mm] \hline \rule{0pt}{3ex}
$3$ & $\geq 3$ & -- \\[2mm] \hline \rule{0pt}{3ex}
$4$ & $2$ &
$b(z^0_1,z^0_2,z^0_3,z^0_4)=b(z^1_1,z^1_2,z^1_3,z^1_4)$
\\[2mm]  \hline \rule{0pt}{3ex}
$4$ & $3$ &
$b_\bbc(z^0_1,z^0_2,z^0_3,z^0_4)=b_\bbc(z^1_1,z^1_2,z^1_3,z^1_4)$
\\[2mm] \hline \rule{0pt}{4ex}
$\geq 4$ & $2$ &
\begin{minipage}[c]{80mm}
\begin{center} \small\sl
$\calo(z^0_1,z^0_2,z^0_3)=\calo(z^1_1,z^1_2,z^1_3)$  et\\[1mm]
%ordre cyclique sur $S^1$ induit par $(z^0_1,z^0_2,z^0_3)$ \'egal\\
%\`a celui induit par $(z^1_1,z^1_2,z^1_3)$ et\\
$b(z^0_1,z^0_2,z^0_3,z^0_j)=b(z^1_1,z^1_2,z^1_3,z^1_k)$
pour $k=4\dots,m$.
\end{center}
\end{minipage} \\[4mm]\hline \rule{0pt}{3ex}
$\geq 4$ & $3$ &
$b_\bbc(z^0_1,z^0_2,z^0_3,z^0_j)=b_\bbc(z^1_1,z^1_2,z^1_3,z^1_k)$
pour $k=4\dots,m$.
\end{tabular}
\end{center}

Les preuves des \'enonc\'es du tableau ci-dessus sont classiques. Pour
$d=2$, on choisit deux transformations de M\"obius
$\sigma_j$ ($j=0,1$) de $\bbr^2\cup\{\infty\}$ telles que
$\sigma_j(z^j_1)=\infty$, $\sigma_j(z^j_2)=0$ et
$\sigma_j(z^j_3)=1$.
Si $\calo(z^0_1,z^0_2,z^0_3)=\calo(z^1_1,z^1_2,z^1_3)$,
alors $\sigma_1^\mun\pcirc\sigma_0\in\moeb{2}$.
De plus, on a que
$$b(z^j_1,z^j_2,z^j_3,z^j_k)=b(\sigma_j(z^j_1),\sigma_j(z^j_2),
\sigma_j(z^j_3),\sigma_j(z^j_k)) = z^j_k,$$
donc l'\'egalit\'e des birapports implique que
$\sigma_1^\mun\pcirc\sigma_0(z^0_k)=z^1_k$ pour $k=4\dots,m$.
La d\'emonstration pour $k=3$ est semblable.
\end{ccote}

\section{Holonomie : le cas g\'en\'erique}\label{Phologen}
Soit $a=(a_1,\dots ,a_m)$ et $b\in\bbr^d$ une valeur r\'eguli\`ere de 
$f^a:(S^{d-1})^m\to\bbr^d$. On d\'esignera par $M_b$
la pr\'eimage $f_a^\mun(b)$ de $b$, qui est donc une
sous-vari\'et\'e de codimension $d$ dans $(S^{d-1})^m$.

Rappelons que les points critiques
de $f^a$ sont les configurations align\'ees, c'est-\`a-dire
les $m$-uples $(z_1,\dots,z_m)$ tels que $z_i=\pm z_j$
\cite[th.\,3.1]{Ha}.
Les valeurs r\'eguli\`eres forment donc un ouvert de $\bbr^d$, 
compl\'ement d'un ensemble fini de sph\`eres centr\'ees en l'origine.
Soit $U$ un ouvert form\'e de valeurs r\'eguli\`eres qui contient $b$.
Soit $\lacets{b}(U)$ l'ensemble des lacets en $b$ dans $U$ qui sont
des courbes lisses par morceau.

\begin{Lemme}\label{ehres}
Pour tout $c\in\lacets{b}(U)$ et tout $q\in M_b$, il existe
une unique $\Delta^a$-courbe $\tilde c$ dans $(S^{d-1})^m$
telle que $f^a\pcirc\tilde c = c$.
\end{Lemme}

\preu
Soit $V:=(f^a)^\mun(U)$. La restriction de $f^a$ \`a $V$ est
donc une submersion qui est propre (puisque $V$ est 
relativement compact). La m\^eme d\'emonstration
que pour une connexion sur un fibr\'e diff\'erentiable
(voir \cite[lemme 2 p. 212]{DNF})
prouve alors le lemme \ref{ehres}. \cqfd

Soit $c$ un lacet \'el\'ement de $\lacets{b}(U)$. A tout $q\in M_b$, 
on peut associer
l'extr\'emit\'e du rel\`evement horizontal de $c$ partant de $q$,
dont l'existence est garantie par le lemme \ref{ehres}.
Ce proc\'ed\'e d\'efinit un
diff\'eomorphisme de $M_b$ appel\'e l'{\it holonomie} de $c$. L'ensemble 
des diff\'eomorphismes obtenus de cette fa\c con constitue un groupe 
$\calh_b(U)$ de diff\'eomomorphismes de $M_b$, appel\'e 
{\it groupe d'holonomie} en $b$ pour l'ouvert $U$. 

La dimension des orbites de l'action de $\calh_b(U)$ sur $M_b$ 
(orbites d'holonomie) d\'epend de $\sharp(a)$. 
Pour un bras articul\'e g\'en\'erique
($\sharp(a)=m$), nous allons d\'emontrer la proposition suivante.

\begin{Proposition}\label{T:transi}
Si $\sharp(a)=m$, le groupe d'holonomie $\calh_b(U)$ 
agit transitivement sur chaque composante connexe de $M_b$.
\end{Proposition}

\preu 
Soit $V:=(f^a)^\mun(U)$ et $q\in V$. 
Soit $\beta_q:G(a,d)\to (S^{d-1})^m$ l'application
$\beta_q(g):=gq$. Si $\sharp(a)=m$, 
l'action du groupe $G(a,d)$ sur $(S^{d-1})^m$ est transitive et
$\beta_q$ est une submersion.

Soit $A\in\lie(G(a,d))$. Le champ de vecteurs
$\alpha_a(A)\in\vect (S^{d-1})^m$ satisfait \`a
$\alpha_a(A)_q=T_{I}\beta_q(A)$ (o\`u $I$
d\'esigne l'\'el\'ement neutre de $G(a,d)$).
Par la proposition \ref{propabd}, l'image 
de $\alpha_a$ est \'egale \`a $\call^a$. On en d\'eduit
que les valeurs en $q$ des champs de $\call^a$
engendrent tout $T_qV$. Par le th\'eor\`eme de
Chow (\cite[\S\,0.4]{Gr}, \cite[th. 2.4 et \S\,2.5]{Be}),
pour tous $q,q'$ dans une m\^eme composante connexe
de $V$, il existe une 
$\Delta^a$-courbe joignant $q$ \`a $q'$ dans $V$. 
Lorsque $q$ et $q'$ sont
dans $M_b$, cette courbe se projette 
sur un lacet dans $\lacets{b}(U)$,
ce qui prouve que $q$ et $q'$ sont dans la m\^eme
orbite d'holonomie. \cqfd

\begin{Remarque}\rm 

a) Lorsque $d=2$, on peut voir plus directement que les
les valeurs en $q$ des champs de $\call^a$
engendrent tout $T_qT^m$ (ou $T_q\bbr^m$). En effet, $\call^a$
contient les champs constants
$A^k$ de \eqref{chpsak}. Si $\sharp(a)=m$, 
ces champs, pour $k=1,\dots ,m$,
sont lin\'eairement ind\'ependants (le calcul aboutit \`a un
d\'eterminant de Vandermonde). 

b) Lorsque $d=2$ la vari\'et\'e $M_b$ peut avoir deux composantes
connexes. Par exemple, si $a:=(\sqrt{3},1)$, $M_b$, pour $b=2$ est
constitu\'e de deux configurations~:
$(e^{i\pi/6},e^{-i\pi/3})$ et $(e^{-i\pi/6},e^{i\pi/3})$.
Comme $G(a,2)$ agit transitivement sur $T^2$, il existe,
par le th\'eor\`eme principal, des
courbes horizontales joignant les
deux points de $M_2$, mais il est clair que
ces courbes passent obligatoirement
par une configuration align\'ee.
Le groupe d'holonomie $\calh_2(U)$ est donc trivial, pour
tout ouvert $U$ contenant $2$ et form\'e de valeurs r\'eguli\`eres.

Plus g\'en\'eralement, soit $a=(a_1,a_2,\varepsilon,\dots,\varepsilon)$
avec $a_1$, $a_2$ et $\varepsilon$ positif, et soit $b\in\bbc$ avec
$|a_1-a_2|<|b|<a_1+a_2$. Si $\varepsilon$ est assez petit, la vari\'et\'e
$M_b$ est diff\'eomorphe \`a deux copies de $T^{m-2}$ et il n'y a pas
de courbe horizontale les joignant sans passer
par une configuration align\'ee.
\end{Remarque}

\section{Holonomie des bras planaires
lorsque $a:=(1,\dots, 1)$}\label{S:orbholo}

Dans tout ce paragraphe, on suppose que $d=2$ et,
sauf indication contraire, $a:=(1,\dots,1)$.
On note $f:=f^a$ et on consid\`ere la pr\'eimage $M_b:=f^\mun(b)$ pour 
une valeur r\'eguli\`ere $b\in\bbc$ de l'application $f$, que l'on voit
donc comme une sous-vari\'et\'e de codimension $2$ dans $T^m$ 
avec la m\'etrique induite. Par abus de notation,
on d\'esignera parfois \'egalement
par $M_b$ la pr\'eimage de $M_b$ dans $\bbr^m$.
Nous allons \'etudier les orbites d'holonomie sur $M_b$.
Soit $\rho : T^m\to S^1$ l'application
$\rho(z_1,\dots ,z_m):=z_1z_2\cdots z_m$. D\'esignons par
$\rho_b:M_b\to S^1$ la restriction de $\rho$ \`a $M_b$.

\begin{Proposition}\label{orbholo}
Soit $U$ un ouvert de $\bbc$ form\'e de valeurs r\'eguli\`eres de $f$
et contenant $b$.
Les orbites de l'action du groupe d'holonomie $\calh_b(U)$ sur $M_b$
sont les trajectoires du champ de gradient de la fonction $\rho_b$. 
\end{Proposition}

\preu  Par le th\'eor\`eme principal, les orbites
d'holonomie sur $M_b$ sont contenues dans l'intersection
avec $M_b$ des orbites de l'action de $G(a,2)=SU(1,1)$ sur $T^m$. 
L'espace tangent en $q$ \`a une telle orbite est
engendr\'e par la valeur en $q$ des champs $C$, $S$ et $A$
de \eqref{chpsak} avec $a_i=1$.
Les vecteurs $C_q$ et $S_q$ sont orthogonaux \`a $T_qM_b$.
Le champ constant $A$ est le gradient de la fonction $\rho$. 
Comme $M_b$ est munie de la m\'etrique induite, la projection
orthogonale de $A_q$ sur $T_qM_b$ donne sur $M_b$ le 
gradient de $\rho_b$. Cela prouve
que les orbites d'holonomie sur $M_b$ sont contenues
dans les trajectoires du gradient de $\rho_b$. 

Pour montrer que les orbites d'holonomie sur $M_b$ sont \'egales
aux trajectoires du gradient de $\rho_b$, on proc\`ede
de la mani\`ere suivante. Pour $s>0$, consid\'erons le chemin $c_s(t)$
dans $\bbc$ ($t\in [0,4s]$) partant de $z$ et parcourant
\`a vitesse $1$ le bord du carr\'e $(z,z+s,z+s+is,z+is)$. 
Si $s$ est assez petit, $c$ est un lacet
\'el\'ement de  $\lacets{b}(U)$.

Le relev\'e horizontal $\tilde c^q(t)$ de $c$ partant de $q\in M_b$
consiste en quatre morceaux de trajectoires des champs
$\tilde e_1,\tilde e_2$, relev\'es horizontaux des champs constants
$e_1=\partial/\partial x$, $e_2=\partial/\partial y\in \vect\,\bbc$.
On a $\tilde c(4s)\in M_b$ et, dans $\bbr^m$, on peut \'ecrire 
\begin{equation}\label{deuord}
\tilde c(4s) = q + \frac{s^2}{2}[\tilde e_1,\tilde e_2]_q + \pcirc (s^2).
\end{equation}
Pour un bras de type $a=(a_1,\dots,a_m)$
d\'ecomposons $[\tilde e_1,\tilde e_2]$ par rapport \`a $C$, $S$ et $A^2$~:
$$[\tilde e_1,\tilde e_2]_q=\alpha_q C_q + \beta_q S_q + \gamma_q A^2_q.$$
Vu l'\'equation \eqref{deuord}, les orbites d'holonomie sur $M_b$ seront
\'egales aux les trajectoires du gradient de $\rho_b$ si 
$\gamma_q\not=0$. Soit $P_q$ la ($2\times 2$)-matrice de Gram
de $D_qf$:
$$P_q:=D_qf(D_qf)^T.$$
On a $\det P_q\not=0$ puisque $q$ est un point r\'egulier de $f$. 
La proposition \ref{orbholo} d\'ecoulera du lemme suivant.

\begin{Lemme}\label{gammaq}
Pour tout bras articul\'e de type $a=(a_1,\dots,a_m)$, on a \\
$\displaystyle \gamma_q=\frac{1}{\det P_q}$.
\end{Lemme}

\preu Soit $f:=f^a$ et soit 
%$$\nabla_q F :=  D_qf^T =(\nabla_q X^a,\nabla_q X^b) = (-AS,AC)$$
%la ``matrice gradient" de $f$. Soit
$V$ un champ de vecteurs sur $\bbc$.
Le champ horizontal $\tilde V$ sur $\bbr^m$ au dessus de $V$ est
donn\'e par 
$$\tilde V_q = (D_qf)^T\, P_q^{-1} V .$$
En effet, on a bien que $\tilde V_q\in\Delta^a_q$ et
$$Df(\tilde V) = Df\, Df^T\, P^{-1} V = P\, P^{-1} V = V .$$
Les cas particuliers $V=e_1=(1,0)^T$ et $V=e_2=(0,1)^T$
donnent
\begin{equation}\label{E1E2} 
\tilde e_1= -Q_{11} AS + Q_{21} AC  \quad , \quad
\tilde e_2= -Q_{12} AS + Q_{22} AC
\end{equation} 
o\`u les $Q_{ij}$ sont les coefficients de la matrice $Q:=P^{-1}$. 
On a donc 
$$[\tilde e_1,\tilde e_2]=\alpha C_q + \beta S_q + \det Q A^2$$
ce qui prouve le lemme \ref{gammaq}. \cqfd

\begin{Remarque}\label{vcrit}\rm
 Un point $q\in\bbr^m$ est r\'egulier si et seulement si 
$\det P_q\not = 0$. Or,
$$ P=\pmatrix{\|AS\|^2 & \langle AS, AC\rangle \cr 
\langle AS, AC\rangle & \|AC\|^2\cr }$$
d'o\`u
$$\begin{array}{ccl} 
\det P_q &=& \displaystyle\sum_{i=1}^n a_i^2\sin^2 q_i \, \sum_{j=1}^n a_j^2\cos^2 q_j 
-\big( \sum_{i=1}^n a_i^2\sin q_i\cos q_i\big)^2 = \\[4 pt] &=&\displaystyle 
\sum_{i\not = j}a_i^2a_j^2 \sin q_i \cos q_j \sin (q_i-q_j) 
= \sum_{i< j}a_i^2a_j^2 \sin^2 (q_i-q_j). 
\end{array}$$
On en d\'eduit que $q$ est un point critique de $f^a$ si et seulement si 
$q_i-q_j=0,\pi$ pour tout $i,j$. Cela signifie que $q$ est une configuration
align\'ee. On retrouve la caract\'erisation des points critiques de
\cite[th.\,3.1]{Ha}.
%En cons\'equence, les valeurs critiques forment un
%ensemble fini de cercles centr\'es en $0$ et les valeurs r\'eguli\`eres
%constituent un ouvert de $\bbc$.
\end{Remarque}

\begin{Proposition}\label{fixholo}
Pour un bras planaire de type $a(1,\dots,1)$,
les points fixes de l'action
du groupe d'holonomie $\calh_b(U)$ sont les configurations
$z=(z_1,\dots ,z_m)\in M_b$ telles que la suite $z_j$ ($j=1,\dots m$) prenne 
exactement 2 valeurs. \end{Proposition}

\preu
Si l'on travaille dans $\bbr^m$,
ce sont les points $q$ tels que $S_q$, $C_q$ et $U_q$ sont
lin\'eairement ind\'ependants, c'est-\`a-dire que la ($m\times 3$)-matrice
construite avec ces trois vecteurs colonne est de rang $\leq 2$.
Cela signifie que pour tout $i,j,k$, on a
$$\left|\matrix{\sin q_i &\cos q_i & 1\cr
\sin q_j &\cos q_j & 1\cr
\sin q_k &\cos q_k & 1\cr}\right| =
\sin (q_i-q_j) + \sin (q_j-q_k) +\sin (q_k-q_i)=0.
$$
Cette \'equation est \'equivalente \`a 
$$\sin \big((q_j-q_k)+(q_k-q_i)\big) = \sin (q_j-q_k) +\sin (q_k-q_i),$$
de la forme $\sin(u+v)=\sin u+\sin v$. Cette derni\`ere \'equation
\'equivaut \`a
\begin{equation}\label{somsin2}
\sin\frac{u+v}{2}\cos\frac{u+v}{2}=\sin\frac{u+v}{2}\cos\frac{u-v}{2}.
\end{equation}
L'\'equation \eqref{somsin2} a lieu si et seulement si
l'une au moins des deux conditions suivantes est v\'erifi\'ee~:
\begin{enumerate}
\item $\cos\frac{u+v}{2}=\cos\frac{u-v}{2}$. Cette condition
est \'equivalente \`a\\ $u\equiv  0\ ({\rm mod\,} 2\pi)$
ou $v\equiv 0\ ({\rm mod\,} 2\pi)$,
c'est \`a dire $e^{iq_i}=e^{iq_k}$ ou $e^{iq_j}=e^{iq_k}$.
\item $\sin\frac{u+v}{2}=0$, c'est-\`a-dire 
$u+v\equiv 0\ ({\rm mod\,} 2\pi)$. Ceci est \'equivalent \`a 
$e^{iq_i}=e^{iq_j}$.
\end{enumerate}
On a ainsi montr\'e que les points fixes du groupe d'holonomie
sont les $q$ tels que, pour tout $i,j,k$, 
l'ensemble $\{e^{iq_i},e^{iq_j},e^{iq_k}\}$
contient au plus deux points. L'ensemble de tous les $e^{iq_j}$ 
doit contenir au moins 2 points,  
sinon on aurait $|b|=m$ et $b$ ne serait pas
une valeur r\'eguli\`ere de $f$. \cqfd

La proposition \ref{fixholo} implique qu'un point critique
$z=(z_1,\dots ,z_m)$ de $\rho_b$ est caract\'eris\'e
par l'ensemble 
$$K=K(z):=\{j\mid Im\langle z_j,b\rangle < 0\}\subset \{1,\dots ,m\},$$
o\`u $\langle,\rangle$ d\'esigne le produit scalaire standard de $\bbc$.
Pour $K$ un sous-ensemble propre de $\{1,\dots ,m\}$,
on d\'esigne par $z^K$ le point critique de $M_b$ tel que
$K(z^K)=K$. Un tel point $z^K$ 
existe dans $M_b$ si et seulement si 
$|m-2|K|\,|<|b|<m$. Par exemple, si $m=4$, le point critique
$z^{\{4\}}$ existe dans $M_b$ si et seulement si $2<|b|<4$.

\begin{Proposition}\label{rhomorse}
La fonction $\rho_b$ est une fonction de
Morse. L'indice du point critique $z^K$ est \'egal \`a $|K|-1$. 
\end{Proposition}

\preu 
Observons que la fonction $\rho:T^m\to S^1$ est satisfait \`a
la condition d'\'equivariance
$$\rho(e^{i\theta}z_1,\dots ,e^{i\theta}z_m)=
e^{i\theta}\rho(z_1,\dots ,z_m).$$
On peut donc,
sans restreindre la g\'en\'eralit\'e, supposer que
$b$ est r\'eel $\geq 0$.
Si $b=0$, il n'y a pas de point critique puisque un tel point
serait alors une configuration align\'ee (le cas $b=0$
sera trait\'e dans l'exemple \ref{E:zero}).
On supposera donc que $b\in\bbr_{>0}$. 

On travaille dans $\bbr^m$ en un point $q^K$ au dessus de $z^K$.
Pour all\'eger la notation, d\'efinissons
$$N:=|K|=\sharp\{j\mid \sin q_j^K<0\} \ \hbox{ et } \ 
P:=\sharp\{j\mid \sin q_j^K>0\}=m-N.$$ 

Observons encore que $M_b$ et $\rho_b$ sont invariantes par l'action
du groupe sym\'etrique sur $T^m$ permutant les coordonn\'ees. Cette action
permute les points critiques ayant des propri\'et\'es analogues. 
On pourra supposer que $K=\{P+1,\dots ,m\}$. 

Au voisinage de $q^K$, la fonction $\rho_b$ admet un rel\`evement 
$\tilde\rho_b$ \`a valeurs 
dans  $\bbr$ de la forme 
\begin{equation}\label{eq:forme}
\tilde\rho_b(q)={\rm cte\,}+\sum_{j=1}^m q_j.
\end{equation}
Notons $\hat f(q)\in\bbc$ l'extr\'emit\'e des $P$ premiers
segments~:
$$\hat f(q):=\sum_{j=1}^{P} e^{iq_j}$$
D\'efinissons $u:=u(q)$, $\xi:=\xi(q)\in ]0,\pi[$, $v:=v(q)$, 
$\eta:=\eta(q)\in ]0,\pi[$
par les formules
$$\hat f(q)= u(q)e^{i\xi(q)}\quad , \quad
b-\hat f(q)= v(q)e^{i(\pi-\eta(q))}.$$
Observons que $\xi$ et $\eta$ sont d\'etermin\'es par $u$ et $v$. En effet,
les points $0,\hat f(q)$ et $b$ forment un triangle de c\^ot\'es $u$, $v$ et $b$
donc, par le th\'eor\`eme du cosinus~:
\begin{equation}\label{cosinus}
\xi(u,v)={\rm arc cos\,}\frac{b^2+u^2-v^2}{2bu} \quad ,\quad
\eta(u,v)={\rm arc cos\,}\frac{b^2-u^2+v^2}{2bv}.
\end{equation}
Consid\'erons 
$(\tilde q,\tilde r)\in\bbr^{P-1}\times\bbr^{N-1}$ d\'efini par
$$\begin{array}{llll}
\tilde q_j&:=&
q_j-\xi(u(q),v(q)) &\hbox{si $j=1,\dots ,\,P-1$}.\\
\tilde r_j&:=&
\eta(u(q),v(q)) + q_{m+1-j} &\hbox{si $j=1,\dots,\,N-1$}.
\end{array}$$
Ces diff\'erentes grandeurs $q_i$, $\tilde q_i$, etc, sont dessin\'ees dans
la figure 1 ci-dessous, dans le cas $m=7$, $P=4$ et $N=3$.

\setlength{\unitlength}{.8mm} 
\begin{picture}(80,30)(8,15)
%\color{gristq}\graphpaper[2](0,-10)(180,100)\color{black}
\put(-3,0){\line(1,0){162}}
\put(5,25){$m=7,\ P=4,\ N=3$}
\savebox{\trait}{\rule{20mm}{0.5mm}}
\savebox{\brasa}{
\put(0,0){\circle*{1.3}}
\put(0,0){\rotatebox{25}{\usebox{\trait}}}
\put(22.6,10.5){\rotatebox{-20}{\usebox{\trait}}}
\put(46.1,2){\rotatebox{20}{\usebox{\trait}}}
\put(69.6,10.5){\rotatebox{-25}{\usebox{\trait}}}
\multiput(0,0)(6,0){16}{\line(1,0){3}}
\put(23,10.6){\line(1,0){18}}
\put(23,10.5){\circle*{1.3}}
\put(46.1,2){\line(1,0){12}}
\put(46.2,2.2){\circle*{1.3}}
\put(70,10.5){\circle*{1.3}}
\put(93,0){\circle*{2.5}}
\qbezier(33,10.6)(34,7)(32,7)
\qbezier(7,0)(8,1.5)(6,3)
\qbezier(55,2.2)(55,2)(54.2,5)
\put(67,1){$u$}
\put(8,1.5){$\tilde q_1$}
\put(34,8){$\tilde q_2$}
\put(56,3){$\tilde q_3$}
}
\put(0,0){\rotatebox{15}{\usebox{\brasa}}}
\savebox{\brasb}{
\put(21,1){$v$}
\put(0,0){\rotatebox{-25}{\usebox{\trait}}}
\put(22.6,-10.5){\rotatebox{44}{\usebox{\trait}}}
\put(22.6,-10.5){\circle*{1.3}}
\put(40.5,7){\rotatebox{-15.4}{\usebox{\trait}}}
\put(40.5,7){\circle*{1.3}}
\put(40.5,7){\line(-1,0){7}}
\put(65,0){\circle*{1.3}}
\qbezier(55,0)(54,3)(56,3)
\qbezier(35,7)(35,4)(36,3)
\multiput(0,0)(6,0){11}{\line(1,0){3}} 
 \put(50.5,1){$\tilde r_1$}
\put(31.5,3.5){$ \tilde r_2$}
}
\put(1,-4){$0$}\put(152,-5){$b$}
\qbezier(16,0)(16,5)(13,9)
\put(16,6.5){$q_1$}
\put(89,24){\rotatebox{-22}{\usebox{\brasb}}}
\put(90,28){$\hat f(q)$}
\put(20,16){\line(1,0){20}}
\qbezier(38,16)(38,15.5)(37.7,15)
\put(42.4,19.4){\vector(-1,-1){4}}
\put(42.5,20){$q_2$}
%\put(45,14){\line(1,0){13}}
%\qbezier(60,14)(60,15.5)(58,18)
%\put(59,17){$q_3$}
\put(134,15.5){\line(1,0){12}}
\qbezier(140,14.5)(140,13)(138,12)
\put(141,12){$ q_7$}
\put(112,6){\line(1,0){12}}
\qbezier(119,6)(119,7)(118,8)
\put(120,7.5){$ q_6$}
\qbezier(24,0)(24,2.8)(22,5)
\put(25,2.4){$\xi$}
\qbezier(136,0)(136,3.5)(138,5)
\put(133.5,2.9){$\eta$}
\put(70,-10){\scshape Figure 1}
\end{picture}
\vskip 25mm

Nous allons param\'etriser $M_b$ au voisinage de $q^K$ 
par des \'el\'ements $(\tilde q,\tilde r)$ situ\'e au voisinage de $0\in
\bbr^{P-1}\times\bbr^{N-1}$.
Observons tout d'abord que $(\tilde q,\tilde r)$ d\'etermine
$u$ et $v$~:
\begin{equation}\begin{array}{lcl}
u(\tilde q)&=&\sum_{j=1}^{P-1}\cos\tilde q_j +
\sqrt{1-(\sum_{j=1}^{P-1}\sin\tilde q_j)^2}  \\[4pt]
v(\tilde r)&=&\sum_{j=1}^{N-1}\cos\tilde r_j +
\sqrt{1-(\sum_{j=1}^{N-1}\sin\tilde r_j)^2}  .
\end{array}\end{equation}
La param\'etrisation $q(\tilde q,\tilde r)$ de $M_b$ au voisinage de 
$q^K=q(0,0)$ est donn\'ee par
\begin{equation}\label{eq:param}
q_i(\tilde q,\tilde r):=\left\{\begin{array}{llllllllll}
\xi(u(\tilde q),v(\tilde r)) + \tilde q_i & \hbox{ si } &
i=1,\dots ,\,P-1\\[2pt]
\xi(u(\tilde q),v(\tilde r)) - 
\arcsin\big(\sum_{j=1}^{P-1}\sin\tilde q_j\big) & \hbox{ si }&
i=P\\[2pt]
-\eta(u(\tilde q),v(\tilde r)) - 
\arcsin\big(\sum_{j=1}^{N-1}\sin\tilde r_j\big) & \hbox{ si }&
i=P+1\\[2pt]
-\eta(u(\tilde q),v(\tilde r)) + \tilde r_{m+1-i} & \hbox{ si } &
i=P+2,\dots m
.\end{array}\right.\end{equation}
On v\'erifie que la matrice jacobienne de l'application $(\tilde q,\tilde
r)\mapsto (q,r)$ est bien de rang $n-2$ (si on lui enl\`eve les lignes $P$ et
$P+1$, on obtient la matrice identit\'e). 

L'expression de $\tilde\rho_b(\tilde q,\tilde r)$ s'obtient des
formules \eqref{eq:forme} et \eqref{eq:param}~:
$$\begin{array}{lclll}
\tilde\rho_b(\tilde q,\tilde r) &=& {\rm cte\,}+
\sum_{j=1}^{P-1}\tilde q_i - 
\arcsin\big(\sum_{j=1}^{P-1}\sin\tilde q_j\big) \quad + \\ &&
\rule{0pt}{13pt}\sum_{j=1}^{N-1}\tilde r_{m+1-i}
-\arcsin\big(\sum_{j=1}^{N-1}\sin\tilde r_j\big) &+&\\ &&
\rule{0pt}{13pt}P\,\xi(u(\tilde q),v(\tilde r))  
-N\,\eta(u(\tilde q),v(\tilde r)) .
\end{array}$$
Comme $\partial u/\partial\tilde q_i (0) =
\partial v/\partial\tilde r_i (0) = 0$,
on a bien $\nabla\tilde\rho_b(0,0)=0$. Quant \`a la matrice hessienne
$\calh\tilde\rho_b$ en $(0,0)$, le calcul de ses coefficients donne
 $$%\begin{equation}
\begin{array}{llll}
\displaystyle\frac{\partial^2 \tilde\rho_b}{\partial\tilde q_k\partial\tilde q_i}(0,0) &=&
s_u\,
\displaystyle\frac{\partial^2 u}{\partial\tilde q_k\partial\tilde q_i}(0,0)
\\[6pt]\rule{0mm}{8mm}
\displaystyle\frac{\partial^2 \tilde\rho_b}{\partial\tilde r_k\partial\tilde r_i}(0,0) &=&
s_v\,
\displaystyle\frac{\partial^2 v}{\partial\tilde r_k\partial\tilde r_i}(0,0)
\\[6pt]\rule{0mm}{8mm}
\displaystyle\frac{\partial^2 \tilde\rho_b}
{\partial\tilde r_k\partial\tilde q_i}(0,0) &=& 0,
\end{array}$$%\end{equation}
avec
$$\begin{array}{rclllllll}
s_u &:=& \displaystyle\frac{\partial}{\partial u}[P\xi
-N \eta](u(0),v(0)) \\[6mm]
s_v &:=& \displaystyle\frac{\partial}{\partial v}[P\xi
-N \eta](u(0),v(0)) \\[6mm]
\displaystyle\frac{\partial^2 u}{\partial\tilde q_k\partial\tilde q_i}(0,0) &=&
\displaystyle\frac{\partial^2 v}{\partial\tilde r_k\partial\tilde r_i}(0,0)
\ =\ -1-\delta_{ik},\end{array}$$
o\`u $ \delta_{ik}$ est le symbole de Kronecker.
Comme $u(0)=P$ et $v(0)=N$, le calcul avec les formules \eqref{cosinus}
donne~:
\begin{equation}\label{eq:finale}\begin{array}{llcllllll}
%\displaystyle\frac{\partial}{\partial u}[P\xi
%-N \eta](u(0),v(0))
s_u &=&
-\displaystyle\frac{PN(1+\cos\gamma)}{2A} & <0 \quad
& i=1\dots P\! -\! 1
\\[6pt]\rule{0mm}{8mm}
%\displaystyle\frac{\partial}{\partial v}[P\xi
%-N \eta](u(0),v(0))
s_v &=&
\displaystyle\frac{PN(1-\cos\gamma)}{2A} & >0
& i=1\dots N\! -\! 1,
\end{array}\end{equation}
\begin{minipage}{9cm}
o\`u $\gamma$ est l'angle en $\hat f(q^K)$ du triangle
$(0,b,\hat f (q^K))$ et $A$ est l'aire de ce triangle; ces formules
s'obtiennent \`a l'aide du th\'eor\`eme du cosinus et de la formule de
H\'eron~: $\scriptstyle 4A=\sqrt{-b^4-P^4-N^4+2b^2P^2+2b^2N^2+2P^2N^2}$. Comme
$b$ est une valeur r\'eguli\`ere, $\gamma\not=0,\pi$, $A\not=0$ et les in\'egalit\'es 
dans \eqref{eq:finale} sont strictes.
\end{minipage}
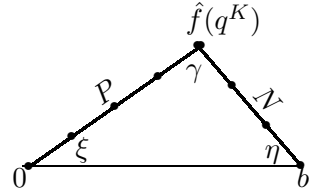
\begin{minipage}{5cm}
\setlength{\unitlength}{.7mm} 
\begin{picture}(60,30)(-8,-5)
%\color{gristq}\graphpaper[2](0,-10)(180,100)\color{black}
\savebox{\brasa}{ 
\thicklines%\linethickness{2pt} 
\put(0,0){\line(1,0){40}} 
\multiput(0,0)(10,0){5}{\circle*{1.5}} 
\put(18,1.5){$P$} 
} 
\savebox{\brasb}{ 
\thicklines%\linethickness{2pt} 
\put(0,0){\line(1,0){30}} 
\multiput(0,0)(10,0){4}{\circle*{1.5}} 
\put(14,1.5){$N$} 
} 
\put(0,0){\rotatebox{35}{\usebox\brasa}} 
\put(35,23){\rotatebox{-49.2}{\usebox\brasb}} 
\thinlines 
\put(2,0){\line(1,0){53}} 
\put(0,-4){$0$}\put(54,-4){$b$} 
%\qbezier(10,0)(10,3)(7,4) 
%\qbezier(49,0)(49,3)(51,6) 
%\qbezier(6,0)(6,3)(5,4)
\put(12,1.7){$\xi$}\put(48,1.5){$\eta$}\put(33,17){$\gamma$}
\put(33,26){$\hat f(q^K)$}
\end{picture}
\end{minipage}

La matrice hessienne $\calh\tilde\rho_b(0,0)$ est donc form\'ee de
deux blocs~:
$$\calh\tilde\rho_b(0,0)\, =\,
\pmatrix{s_uM(N) & 0\cr 0 & s_vM(N)\cr},$$
o\`u $M(k)$ est la $(k\times k)$-matrice dont les coefficients diagonaux sont
$-2$ et les autres $-1$. Comme $M(k)$ est somme de la matrice $-I$
et d'une matrice de rang $1$, ses valeurs propres sont faciles \`a calculer~:
l'une vaut $-(k+1)$ et toutes les autres valent $-1$. En particulier,
elles sont toutes n\'egatives. On d\'eduit ainsi de
\eqref{eq:finale} que
$\calh\tilde\rho_b(0,0)$ a $N-1$ valeurs propres n\'egatives
et $P-1$ valeurs propres positives. Le point $z^K$ est donc
un point critique non-d\'eg\'en\'er\'e de $\rho_b$ d'indice $|K|-1$.
\cqfd

\begin{Remarque}\label{R:niveau} Tous les points critiques 
d'un m\^eme indice sont sur un m\^eme niveau de $\rho_b$.
\end{Remarque}

\paragraph{Exemples}{ \bf :}

\begin{ccote}\label{E:sph} \rm
Pour le bras articul\'e avec $m$ segments tous de longueur $1$,
la vari\'et\'e $M_b$ est diff\'eomorphe \`a l'espace de polygones
$\nua{m+1}{2}(1,\dots,|b|)$
(voir \cite{HR}). Lorsque $m-2<|b|<m$, vari\'et\'e $M_b$ est donc
diff\'eomorphe \`a la sph\`ere $S^{m-2}$ \cite[Exemple 6.5]{HR}.
Les angles $q_i$ sont bien d\'etermin\'es 
dans l'intervalle ouvert $(-\pi,\pi)$ (aucun segment ne pouvant se
``retourner'') et on pourra donc 
identifier $M_b$ \`a la 
composante connexe de $\tilde M_b\subset\bbr^m$ 
pour laquelle $|q_j|<\pi,\ \forall\! j$.
De m\^eme, l'application 
$\rho_b :M_b\to S^1$ se rel\`eve en $\tilde\rho : M_b\to\bbr$ d\'efinie
par $\tilde\rho_b(q)=\sum_{i=1}^mq_i$ (un tel rel\`evement continu n'existe
que si $|b|>m-2$).

Les points critiques $z^K$ existent dans $M_b$ pour tout sous-ensemble
propre $K$ de $\{1,\dots ,m\}$. 
Par le th\'eor\`eme \ref{rhomorse},
pour tout $1\leq k\leq m-1$, la fonction de Morse $\tilde\rho_b$ poss\`ede
donc $\big(\,{}^m_k\big)$ points critiques d'indice $k-1$. 

Pour visualiser le champ $\grad\,\tilde\rho_b= \grad\,\rho_b$
sur $M_b\iso S^{m-2}$, 
on consid\`ere une triangulation
$\calt$ sur $S^{m-2}$ isomorphe \`a la premi\`ere subdivision barycentrique 
du bord du $(m-1)$-simplexe $\Delta^{m-1}$. 
Si l'on num\'erote les sommets de $\Delta^{m-1}$ de $1$ \`a $m$, les 
sommets de $\calt$ sont index\'es par les sous-ensembles propres $K$
de $\{1,\dots m\}$. 
On peut v\'erifier que $\grad\rho_b$ 
est le champ associ\'e \`a $\calt$
(voir \cite[p. 611--612]{Sp})~:
le barycentre d'un $k$-simplexe est un z\'ero de $\grad\rho_b$
correspondant \`a un point critique d'indice $k$ de $\rho_b$.
\end{ccote}

\noindent\begin{minipage}{6.5cm}
\begin{ccote}\rm Consid\'erons le cas particulier de l'exemple
\ref{E:sph} o\`u $m=3$. Pour $1<b<3$,
la vari\'et\'e $M_b$ est diff\'eomorphe \`a un cercle et
la fonction $\rho_b$ admet 3 points critiques d'indice $0$
et 3 d'indice $1$. La figure ci-contre montre les configurations
critiques. Les fl\`eches indiquent le sens des trajectoires de
$\grad\rho_b$, allant des points critiques d'indice $0$
\`a ceux d'indice $1$.
\end{ccote}
\end{minipage}
\begin{minipage}{6cm}
\begin{pspicture}(0,0)(0,0)
\setlength{\unitlength}{8mm}
\put(5,0){\pscircle{2}}
\put(7,0){\circle*{.2}}
\put(3,0){\circle*{.2}}
\put(6,1.73){\circle*{.2}}
\put(6,-1.73){\circle*{.2}}
\put(4,1.73){\circle*{.2}}
\put(4,-1.73){\circle*{.2}}
%box1
\put(7.5,-0.2){\rotatebox{-18}{
\put(0,0){\line(1,1){.5}}
\put(0.5,0.5){\line(1,1){.5}}
\put(1,1){\line(1,-1){.5}}
\put(0,0){\circle*{.1}}
\put(0.5,0.5){\circle*{.1}}
\put(1,1){\circle*{.1}}
\put(1.5,0.5){\circle*{.1}}
}}
%box4
\put(1,0.2){\rotatebox{17}{
\put(0,0){\line(1,-1){.5}}
\put(0.5,-0.5){\line(1,-1){.5}}
\put(1,-1){\line(1,1){.5}}
\put(0,0){\circle*{.1}}
\put(0.5,-0.5){\circle*{.1}}
\put(1,-1){\circle*{.1}}
\put(1.5,-0.5){\circle*{.1}}
}}
%box2
\put(6.5,1.7){\rotatebox{-71}{
\put(0,0){\line(-1,1){.5}}
\put(-0.5,0.5){\line(1,1){.5}}
\put(0,1){\line(1,1){.5}}
\put(0,0){\circle*{.1}}
\put(-0.5,0.5){\circle*{.1}}
\put(0,1){\circle*{.1}}
\put(0.5,1.5){\circle*{.1}}
}}
%box6
\put(6.3,-1.9){\rotatebox{16}{
\put(0,0){\line(1,-1){.5}}
\put(0.5,-0.5){\line(1,1){.5}}
\put(1,0){\line(1,-1){.5}}
\put(0,0){\circle*{.1}}
\put(0.5,-0.5){\circle*{.1}}
\put(1,0){\circle*{.1}}
\put(1.5,-0.5){\circle*{.1}}
}}
%box3
\put(1.9,1.9){\rotatebox{-16}{
\put(0,0){\line(1,1){.5}}
\put(0.5,0.5){\line(1,-1){.5}}
\put(1,0){\line(1,1){.5}}
\put(0,0){\circle*{.1}}
\put(0.5,0.5){\circle*{.1}}
\put(1,0){\circle*{.1}}
\put(1.5,0.5){\circle*{.1}}
}}
%box5
\put(1.7,-1.8){\rotatebox{-18}{
\put(0,0){\line(1,-1){.5}}
\put(0.5,-0.5){\line(1,1){.5}}
\put(1,0){\line(1,1){.5}}
\put(0,0){\circle*{.1}}
\put(0.5,-0.5){\circle*{.1}}
\put(1,0){\circle*{.1}}
\put(1.5,0.5){\circle*{.1}}
}}
\thicklines
\put(6.735,.9){\vector(-1,3){.1}}
\put(6.735,-.9){\vector(-1,-3){.1}}
\put(3.24,.9){\vector(-1,-3){.1}}
\put(3.24,-.9){\vector(-1,3){.1}}
\put(5,2){\vector(1,0){.1}}
\put(5,-2){\vector(1,0){.1}}

\put(6,-0.2){$z^{\{3\}}$}
\put(5.1,1){$z^{\{2,3\}}$}
\put(3.6,1.1){$z^{\{2\}}$}
\put(3.2,-0.2){$z^{\{1,2\}}$}
\put(3.7,-1.4){$z^{\{1\}}$}
\put(5.1,-1.5){$z^{\{1,3\}}$}
\end{pspicture}
\end{minipage}

\begin{ccote}\rm Consid\'erons le bras articul\'e \`a 4 segments~:
$a=(1,1,1,1)$. Les valeurs critiques de $f_a$ sont $0$ et les cercles
de rayon $2$ et $4$. Si $2<|b|<4$, la vari\'et\'e $M_b$ est diff\'eomorphe
\`a $S^2$ et on est dans le cas de l'exemple \ref{E:sph}. Lorsque 
$0<|b|<2$, la vari\'et\'e $M_b$ est diff\'eomorphe \`a l'espace des configurations
planaires (modulo rotations) d'un pentagone r\'egulier; cet espace est
diff\'eomorphe \`a une surface de genre $4$ (voir \cite[Table V]{HR}).
L'application $\rho_b$ a alors 6 points critiques, 
qui sont tous d'indice 1 (points ``selles'', 
correspondant aux configurations~:
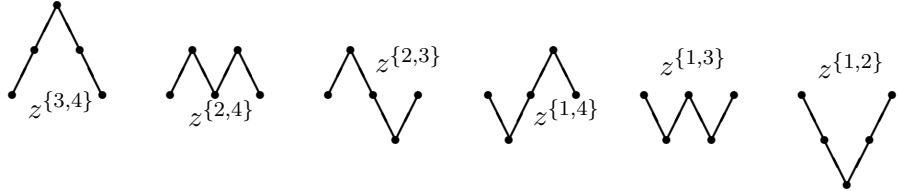
\begin{minipage}{\hsize}
\setlength{\unitlength}{.15mm} 
\begin{picture}(60,30)(-25,80)
\savebox{\brasa}{ 
\thicklines%\linethickness{2pt} 
\put(0,0){\line(1,2){20}} 
\put(0,0){\circle*{6}}} 
\savebox{\brasb}{ 
\thicklines%\linethickness{2pt} 
\put(0,0){\line(1,-2){20}} 
\put(0,0){\circle*{6}}}
\savebox{\conab}{
\put(0,0){\usebox\brasa}
\put(20,40){\usebox\brasa}
\put(40,80){\usebox\brasb}
\put(60,40){\usebox\brasb}
\put(80,0){\circle*{6}}
\put(14,-20){$z^{\{3,4\}}$}}
\savebox{\conac}{
\put(0,0){\usebox\brasa}
\put(20,40){\usebox\brasb}
\put(40,0){\usebox\brasa}
\put(60,40){\usebox\brasb}
\put(80,0){\circle*{6}}
\put(16,-28){$z^{\{2,4\}}$}}
\savebox{\conad}{
\put(0,0){\usebox\brasa}
\put(20,40){\usebox\brasb}
\put(40,0){\usebox\brasb}
\put(60,-40){\usebox\brasa}
\put(80,0){\circle*{6}}
\put(42,20){$z^{\{2,3\}}$}}
\savebox{\conbc}{
\put(2,0){\usebox\brasb}
\put(20,-40){\usebox\brasa}
\put(40,0){\usebox\brasa}
\put(60,40){\usebox\brasb}
\put(80,0){\circle*{6}}
\put(42,-27){$z^{\{1,4\}}$}}
\savebox{\conbd}{
\put(0,0){\usebox\brasb}
\put(20,-40){\usebox\brasa}
\put(40,0){\usebox\brasb}
\put(60,-40){\usebox\brasa}
\put(80,0){\circle*{6}}
\put(14,18){$z^{\{1,3\}}$}}
\savebox{\concd}{
\put(0,0){\usebox\brasb}
\put(20,-40){\usebox\brasb}
\put(40,-80){\usebox\brasa}
\put(60,-40){\usebox\brasa}
\put(80,0){\circle*{6}}
\put(14,15){$z^{\{1,2\}}$}}
\put(0,0){\usebox\conab}
\put(140,0){\usebox\conac}
\put(280,0){\usebox\conad}
\put(420,0){\usebox\conbc}
\put(560,0){\usebox\conbd}
\put(700,0){\usebox\concd}
\end{picture}\end{minipage}
\vskip 25mm
Rappelons (voir remarque \ref{R:niveau})
que ces 6 points critiques sont sur
un m\^eme niveau~: $\rho_b(z^{i,j})=1$ pour tout $i,j$. La visualisation
des trajectoires de $\grad\rho_b$ sur la surface $M_b$ de genre 4
est un int\'eressant d\'efi. 
\end{ccote}

\begin{ccote}\label{E:zero}\rm
Lorsque $m$ est impair, le disque ouvert de rayon 1
centr\'e en $0$ est form\'e de valeurs r\'eguli\`eres.
Si $|b|<1$, il n'y a aucune configuration satisfaisant 
la condition de proposition \ref{fixholo}. La fonction
$\rho_b$ n'a donc aucun point critique et est une submersion.
Le groupe d'holonomie $\calh_b$ agit sans point fixe sur $M_b$.

Dans le cas particulier $z=0$, on a
$\grad\rho_b=U$ et le flot de $\grad\rho_b$ est p\'eriodique,
donn\'e par l'action diagonale de $S^1$ sur $T^m$.
Chaque configuration effectue 
simplement un mouvement de
rotation \`a vitesse angulaire constante. 
Cela prouve que $M_0$, et donc $M_b$ pour $|b|<1$, est 
diff\'eomorphe \`a $M^{(\bar a)}_1\times S^1$, o\`u 
$\bar a$ est le bras articul\'e \`a $m-1$ segments, tous de 
longueur $1$. 
\end{ccote}

\vskip .3 truecm\goodbreak

\vskip .5 truecm\small
\noindent \parbox[t]{6 truecm}{Jean-Claude HAUSMANN\\
Math\'ematiques-Universit\'e\\ B.P. 240, \\
 CH-1211 Gen\`eve
 24, Suisse\\ hausmann@math.unige.ch} \ \hfill \hfill \

\end{document}